\documentclass[11pt,thmsa,a4paper,reqno]{article}
\usepackage{scrextend}
\changefontsizes{12pt}
\usepackage[all]{nowidow}
\usepackage{float}
\usepackage{xfrac}
\usepackage[utf8x]{inputenc}
\usepackage{amscd,amssymb,amsmath}
\usepackage{amsrefs}
\usepackage{tikz-cd}
\usepackage{bigdelim}
\usepackage{logicproof}
\usepackage{longtable}
\usepackage{parskip}
\usepackage{amsthm}
\usepackage{thmtools}
\usepackage[scr,scaled=1.1]{rsfso}

\makeatother

\let\longlogicproof\logicproof
\let\endlonglogicproof\endlogicproof
\patchcmd{\longlogicproof}
  {\center\expandafter\tabular}
  {\expandafter\longtable}
  {}{}
\patchcmd{\endlonglogicproof}
  {\endcenter}
  {}
  {}{}
\patchcmd{\endlonglogicproof}
  {\endtabular}
  {\endlongtable}
  {}{}

\usepackage{mathptmx}
\usepackage[T1]{fontenc}

\usepackage{bm}
\usepackage{bbm}

\theoremstyle{definition}

\theoremstyle{definition}

\theoremstyle{remark}

\usepackage{graphicx}
\usepackage{resizegather}
\usepackage{arydshln}
\usepackage{soul}
\usepackage{relsize,etoolbox}
 \usepackage[all,cmtip,2cell]{xy}
\UseAllTwocells
\xyoption{2cell}

\usepackage{tikz}
\usepackage{extarrows}
\usepackage{amssymb}

\numberwithin{equation}{section}

\usepackage[pdftex,
                paper=a4paper,
                portrait=true,
                textwidth=165mm,
                textheight=235mm,
                tmargin=2.5cm,
                marginratio=1:1]{geometry}

\newtheoremstyle{mydef}
  {}		
  {}		
  {}		
  {}		
  {\scshape}	
  {. }		
  { }		
  {\thmname{#1}\thmnumber{ #2}\thmnote{ #3}}	

\theoremstyle{definition}
\newtheorem{theorem}{Theorem}[section]
\newtheorem*{theorem*}{Theorem}

\newtheorem{proposition}[theorem]{Proposition}
\newtheorem*{proposition*}{Proposition}
\newtheorem{lemma}[theorem]{Lemma}
\newtheorem*{lemma*}{Lemma}
\newtheorem{corollary}[theorem]{Corollary}
\newtheorem*{corollary*}{Corollary}
\theoremstyle{definition}
\newtheorem{definition}[theorem]{Definition}
\newtheorem{example}[theorem]{Example}
\newtheorem{notation}[theorem]{Notation}
\theoremstyle{remark}
\newtheorem{remark}[theorem]{Remark}

\DeclareMathAlphabet{\mathsfit}{T1}{\sfdefault}{\mddefault}{\sldefault}
\SetMathAlphabet{\mathsfit}{bold}{T1}{\sfdefault}{\bfdefault}{\sldefault}

\DeclareFontFamily{U}{matha}{\hyphenchar\font45}
\DeclareFontShape{U}{matha}{m}{n}{
      <5> <6> <7> <8> <9> <10> gen * matha
      <10.95> matha10 <12> <14.4> <17.28> <20.74> <24.88> matha12
      }{}
\DeclareSymbolFont{matha}{U}{matha}{m}{n}
\DeclareFontSubstitution{U}{matha}{m}{n}
\DeclareMathSymbol{\abxcup}{\mathbin}{matha}{'131}

\DeclareMathAlphabet{\mathbfit}{OML}{cmm}{b}{it}

\let\longlogicproof\logicproof
\let\endlonglogicproof\endlogicproof
\patchcmd{\longlogicproof}
  {\center\expandafter\tabular}
  {\expandafter\longtable}
  {}{}
\patchcmd{\endlonglogicproof}
  {\endcenter}
  {}
  {}{}
\patchcmd{\endlonglogicproof}
  {\endtabular}
  {\endlongtable}
  {}{}

\begin{document}

\title{Propositional Calculus with Multiple Negations}
\author{\textsc{Oscar Ramirez} 
\footnote{Email:  odramirezc@gmail.com, Bogotá D.C., Colombia.}}
\date{}
\maketitle

\begin{abstract} 
One advantage of paraconsistent logic is that it can deal with inconsistencies without making the system trivial. However, unlike classical propositional calculus, its deductive system is limited, and the meaning of paraconsistent negation is still not clear. This article presents a logical system that brings together the strengths of both approaches. The Propositional Calculus with Multiple Negations $\left(\textbf{CPN}_{n}\right)$ is a generalization of classical propositional logic in which a finite number of negations (each weaker than the classical one but with similar behavior) are added. This makes it possible to introduce weak inconsistencies in a controlled way without leading to triviality.
\end{abstract}

\maketitle
\setcounter{tocdepth}{1}
\tableofcontents

\section{Introduction}
The propositional calculus with multiple negations, or simply $\textbf{CPN}_{n}$ (where $n$ is the number of negations), addresses the following questions: Can a proposition be true for one observer and false for another? Can truth and falsity be seen as local notions? And if so, how can such an idea be formally represented? While these questions may receive positive answers from a philosophical point of view, their logical and mathematical treatment is much more complex. Classical logic, as well as many-valued logic, allows propositions to take two or more truth values \cite{caicedo}, \cite{Hod}, \cite{Men}; however, this does not explain whether something true under one criterion must also be true under another. In $\textbf{CPN}_{n}$, statements pass through different modalities of truth across various possible worlds. To show how this system works, we present two representative examples.

Let us consider three individuals: Bob, Alice, and Mitchell. Suppose we want to evaluate the proposition $\psi$: “Bob is a good person.” If Bob meets some basic social standards, $\psi$ may be true for Alice but not for Mitchell. What Alice considers good is not necessarily what Mitchell considers good, and vice versa. To give a logical interpretation of this scenario, let us define $\mathcal{M}_{A}= \{T_{A}, F_{A} \}$ as Alice’s set of truth values, and $\mathcal{M}_{M}=\{ T_{M}, F_{M} \}$ as Mitchell’s set of truth values. The truth table for $\psi$ can then be written as follows:

\begin{table}[h]
\label{table: truthtablealicebob}
    \centering

    \begin{tabular}{cl}
        $\psi$ \\
        \cline{1-1}
       
        $T_{A}$ & \rdelim\}{2}{*}[World of Alice]\\
        $F_{A}$ & \\
        \cline{1-1}
        $T_{M}$ & \rdelim\}{2}{*}[World of Mitchell]\\
        $F_{M}$ & \\
    \end{tabular}
    \vspace{0.2cm}
    \caption{Truth table of $\psi$}
\end{table}

According to the truth table, the possible truth values for $\psi$ are $\left( T_{A}, T_{M} \right)$, $\left( T_{A}, F_{M} \right)$, $\left(F_{A}, T_{M} \right)$, or $\left(F_{A}, F_{M} \right)$. The next example shows a case where the law of non-contradiction can fail under weaker forms of negation. Consider the proposition $\varphi$: “The cat is alive.” Let $\mathcal{M}_{1}=\{T_{1}, F_{1} \}$ and $\mathcal{M}_{2}=\{T_{2}, F_{2} \}$ be two possible worlds in which $\varphi$ is evaluated. Because the worlds differ, their negations differ as well; write $\neg_{1}$ and $\neg_{2}$ for the negation operations in $\mathcal{M}_{1}$ and $\mathcal{M}_{2}$, respectively. With this setup, we analyze the formula $\varphi \land \neg_{1}\varphi$ by means of a truth table.

\begin{table}[h]
\label{table: truthtableforcontradiction}
\centering
\begin{tabular}{cc|c}
$\varphi$ & $\neg_{1}\varphi$ & $\varphi \land\neg_{1}\varphi$ \\
\hline
$T_{1}$ & $F_{1}$ & $F_{1}$\\
$F_{1}$ & $T_{1}$ & $F_{1}$ \\
$T_{2}$ & $T_{2}$ & $T_{2}$ \\
$F_{2}$ & $F_{2}$ & $F_{2}$ \\
\end{tabular}
\vspace{0.2cm}
\caption{Truth table for $\varphi \land \neg_{1}\varphi$}
\end{table}

In Table 2, some points remain unclear, such as the interpretation of $\land$ in this system. We will clarify these questions in the course of the article. Without full rigor at this stage, we can already see that in the world $\mathcal{M}_{2}$ it is possible for $\varphi \land \neg_{1}\varphi$ to be true (and therefore for $\neg_{1} \left( \varphi \land \neg_{1}\varphi \right)$ to be false) in $\mathcal{M}_{2}$. Moreover, it may happen that $\neg_{1}\varphi$ is true in $\mathcal{M}_{1}$ while $\varphi$ is true in $\mathcal{M}_{2}$.

As seen above, the system $\textbf{CPN}_{n}$ is a generalization of classical propositional calculus, in which a proposition’s truth value may depend on the number of worlds in which it is evaluated. Thus, a proposition that is true in one world may be false in another. The previous example relates to themes in physics, but this paper does not take a position on the interpretation of quantum mechanics.

\textit{Content of the paper.} In Section 2, we introduce the preliminary notion of chains over the natural numbers. In Section 3, we present the deductive system for $\textbf{CPN}_{n}$, including its axioms, rules of inference, and some relevant results. In Section 4, we give the semantics for $\textbf{CPN}_{n}$, prove the soundness of the axioms, and state the Soundness Theorem. We then show that certain formulas from classical propositional calculus that use classical negation no longer hold in $\textbf{CPN}_{n}$ when classical negation is replaced by weak negations. In the final section, we provide the proof of the Completeness Theorem and describe an interesting relationship among the different logical systems $\textbf{CPN}_{n}$.

\section{Preliminaries}
\subsection*{Chains}\label{chains} 
If a proposition $\varphi$ is evaluated across several states or worlds, we may need to negate it in one or more of them at the same time. For this purpose, we use the notion of a chain whose symbols belong to a subset of the natural numbers. For background, see \cite{Korgi} and \cite{Sipser}.

\begin{definition}
 \label{def:chains}  

A \emph{chain} of length $k$ is a finite sequence of $k$ distinct natural numbers. We write $c_{k}$ to denote any chain.
\end{definition} 

Although the notion of a chain introduced here uses the set of natural numbers as its alphabet, from now on we consider only chains formed from finite subsets of natural numbers, specifically subsets of the form $[n]=\{1, \ldots, n \}$. We also allow the empty chain, that is, the chain with no symbols from the alphabet, denoted by $\epsilon$.

\begin{definition}
    \label{def: setofallchains}
The set of all chains $c_{k}$ over $[n]$, including the empty chain $\epsilon$, is denoted by $[n]^{**}$.     
\end{definition}

\begin{example}
\label{exa: setofchains}
 For $[3]$, we have that: $$[3]^{**} = \{\epsilon,1,2,3,12,13,21,23,31,32,123,132,213,231,321,312 \}\text{.}$$ 
\end{example}

Since chains cannot contain repeated symbols from $[n]$, the set $[n]^{**}$ is finite. However, we will restrict attention to a smaller subset, because we are not interested in all chains that share the same length and the same symbols. For example, when $[3]$ is the alphabet, we do not distinguish between the chains $312$ and $123$. To address this, we define an equivalence relation on $[n]^{**}$ and remove the “repeated” elements.

\begin{proposition}
\label{prop: relationofequivalence}

Let $n\in\mathbb{N}$, and let $[n]^{**}$ be the set of all chains over $[n]$. The relation on $[n]^{**}$ defined by 
$$c_{k} \equiv c_{r} \hspace{0.2cm} \textit{if and only if } k=r \hspace{0.2cm} \textit{and the chains consist of the same symbols }    $$

is an equivalence relation.     
\end{proposition}

\begin{remark}
 \label{rem: relationofequivalence}  

With a slight abuse of notation, we write $c_{k}$ instead of $[c_{k}]$ to refer to the equivalence class. Moreover, the quotient set $\sfrac{[n]^{**}}{\equiv}$ will, from now on, be denoted by $[n]^{*}$. Additionally, the chain of length $n$, that is, the one containing all the symbols of the alphabet $[n]$, will be abbreviated as $(n)$.
\end{remark}

\begin{example}
\label{ref: abusingnotation}
For $[3]$, we have that: $$[3]^{*}=\{ \epsilon,1,2,3,12,13,23,(3) \}\text{.}$$    
\end{example}

\begin{definition}
  \label{def:concatenation}  

Let $c_{k}$ and $c_{r}$ be two chains over $[n]$ and $[m]$, respectively. The \textit{concatenation} of $c_{k}$ and $c_{r}$, denoted $c_{k} \cdot c_{r}$, is defined as follows:

\begin{itemize}
    \item[$(i)$] If $c_{r}=\epsilon$, then $c_{k}\cdot\epsilon= \epsilon \cdot c_{k}=c_{k}$; 
    \item[$(ii)$] $c_{k}\cdot c_{r}=c_{s}$, where $c_{s}$ is the chain obtained by writing the symbols of $c_{k}$ followed by those of $c_{r}$, with the convention that if a symbol occurs in both chains, it appears only once in $c_{s}$.
\end{itemize}

\end{definition}

\begin{definition}
    \label{def: subchains}
A chain $c_{r}$ is a \textit{subchain} of $c_{k}$ if $r\leq k$ and there exist chains $x_{s}$ and $y_{s'}$ such that $c_{k}=x_{s}\cdot c_{r}\cdot y_{s'}$. 
\end{definition}

\begin{definition}
  \label{def:coconcatenation}  

Let $c_{k}$ and $c_{r}$ be two chains over $[n]$ and $[m]$, respectively. The \textit{coconcatenation} of $c_{k}$ and $c_{r}$, denoted $c_{k} \otimes c_{r}$, is defined as follows:

\begin{itemize}
    \item[$(i)$] If $c_{r}=\epsilon$, then $c_{k}\otimes \epsilon= \epsilon \otimes c_{k}=c_{k}$; 
    \item[$(ii)$] $c_{k}\otimes c_{r}=c_{s}$, where $c_{s}$ is the chain formed by writing the symbols of $c_{k}$ followed by those of $c_{r}$ and deleting any symbol that occurs in both chains. Equivalently, $c_{s}$ keeps only the symbols that occur in exactly one of the two chains.
 
\end{itemize}

\end{definition}

\begin{example}
\label{exa: examplecococatenaction}
Consider the chains $c_{3}=123$ over the alphabet $[3]$ and $c_{5}=12456$ over the alphabet $[6]$. Then $c_{3}\otimes c_{5}=3456$.  
\end{example}

\begin{notation}
    \label{not: notationchain}
From now on, the symbols of a chain will be written with commas. For example, the chain $c_{5}=2345$ over $[5]$ will be written as $c_{5}=2,3,4,5$.
\end{notation}

\begin{remark}
\label{rem: conmutativechain}
 The definition \ref{def: subchains} corresponds to the classical notion of a subchain. However, since we are working with equivalence classes, the concatenation and coconcatenation operations are commutative, so the order of the factors $c_{r}$, $x_{s}$, and $y_{s'}$ does not matter.
\end{remark}

When two chains over $[n]$ are disjoint, that is, they share no symbols from the alphabet, concatenation and coconcatenation coincide. Moreover, the coconcatenation of two chains is always a subchain of their concatenation.

\begin{proposition}
\label{prop: chainandcomplementchain}
Let $c_{k}$ be a chain over the alphabet $[n]$, for some $n\in \mathbb{N}$. Then there exists a chain $c_{n-k}^{'}$, with the same alphabet $[n]$, such that

 $$c_{k}\cdot c_{n-k}^{'}=c_{n-k}^{'} \cdot c_{k}=c_{k} \otimes c_{n-k}^{'} = c_{n-k}^{'}\otimes c_{k}=(n)$$

 The chain $c_{n-k}^{'}$ is called the complementary chain of $c_{k}$.
\end{proposition}

\begin{proposition}
\label{prop: propietiesofproductofchains}
Let $c_{k}$ be an arbitrary chain over the alphabet $[n]$, where $n\in\mathbb{N}$. Then we have: 
\begin{enumerate}
    \item[\textup{(i)}] $c_{k}\otimes c_{k}= \epsilon \text{;}$  
    \item[\textup{(ii)}] $(n)\otimes c_{k}=c_{k}\otimes(n) =c_{n-k}^{'} \text{;}$
    \item[\textup{(iii)}] $(n) \otimes c_{n-k}^{'}=c_{n-k}^{'}\otimes (n) =c_{k} \text{.}$ 
   
\end{enumerate}
\end{proposition}

\section{Deductive System of $\textbf{CPN}_{n}$}

The behavior of $\textbf{CPN}_{n}$ is similar to classical propositional calculus. More importantly, $\textbf{CPN}_{n}$ is essentially the same as $\textbf{CPC}$\footnote{Classical propositional calculus.}, except that it introduces weaker forms of negation.

\textbf{Syntax.}\label{syntax} For each natural number $n$, the alphabet of the propositional calculus with $n$ negations consists of the following symbols:

\begin{enumerate}
   \item[1.] An enumerable set $P_{n}$ of atomic formulas $\varphi_{1},\ldots\varphi_{m},\ldots$;
    \item[2.]  A set of symbols $\mathcal{C}_{n}$ where $$\mathcal{C}_{n}=\{ \perp_{c_{k}} : c_{k} \ \text{is a chain over} \  [n]  
 \ \textrm{and} \ 1\leq k \leq n-1  \}\text{;}$$
     \item[3.] A finite set of negations $\{ \neg_{c_{k}} \}_{1\leq k\leq n}$, where each chain $c_{k}$ is over the alphabet $[n]$, and a connective $\longrightarrow_{(n)}$, called the $n$-implication;
     \item[4.] Two symbols $\perp_{(n)}$ and $\top_{(n)}$ representing contradiction and truth, respectively;
     \item[5.] Bracket symbols $\left(,\right)$.
\end{enumerate}

\begin{remark}
\label{rem: observationnnegation}
In item 3, when the length $k$ equals $n$, the negation $\neg_{c_{n}}$ is called the $n$-negation and is denoted by $\neg_{(n)}$. In $\textbf{CPN}_{n}$ the $n$-negation is also called strong negation.
\end{remark}

\begin{definition}
\label{def: wellformedformula}
A string of symbols from $\Omega_{n}$, where $\Omega_{n} = P_{n} \cup \mathcal{C}_{n}$, is a well-formed formula if and only if it is obtained by the following rules:

\begin{enumerate}
    \item The elements of $\Omega_{n}$ are well-formed formulas;
    \item If $\varphi$ is a well-formed formula, then $\neg_{c_{k}}\varphi$ is also a well-formed formula for $1 \leq k \leq n$;
    \item If $\varphi$ and $\psi$ are well-formed formulas, then $\varphi \longrightarrow_{(n)} \psi$ is a well-formed formula.
\end{enumerate}

The set of all well-formed formulas is denoted by $\mathcal{F}(\Omega_{n})$.
\end{definition}

\begin{definition}
\label{def: nconnectives}
For each $n\in\mathbb{N}$ and for all $\varphi, \psi \in \mathcal{F}\left(\Omega_{n} \right)$, the connectives $\land_{(n)}$, $\lor_{(n)}$, and $\longleftrightarrow_{(n)}$ are defined as follows:

\begin{enumerate}
    \item $\varphi \land_{(n)} \psi := \neg_{(n)} \left( \varphi \longrightarrow_{(n)} \psi \right) \text{;}$ 
    \item $\varphi \lor_{(n)} \psi:= \neg_{(n)}\varphi \longrightarrow_{(n)} \psi \text{;}$
    \item $\varphi \longleftrightarrow_{(n)} \psi:= \left( \varphi \longrightarrow_{(n)} \psi \right) \land_{(n)} \left( \psi \longrightarrow_{(n)} \varphi \right) \text{.}$
\end{enumerate}
The connectives $\land_{(n)}$, $\lor_{(n)}$, and $\longleftrightarrow_{(n)}$ are called the $n$-conjunction, $n$-disjunction, and $n$-biconditional, respectively.
\end{definition}

\begin{notation}
\label{not: notationnegationandocntradiction}
For every well-formed formula $\varphi$, we define $\neg_{\epsilon} \varphi:= \varphi$ and $\perp_{\epsilon}:= \top_{(n)}\text{.}$
\end{notation}


\textbf{Axioms.}\label{axioms} For $n\in\mathbb{N}$, an axiom of $\textbf{CPN}_{n}$ is any formula of the following forms, where $\varphi$, $\psi$, and $\chi$ are arbitrary formulas in $\mathcal{F}\left( \Omega_{n} \right)$:

\begin{tabbing}
  \textbf{A1}\quad \= $\varphi \longrightarrow_{(n)} \left( \psi \longrightarrow_{(n)} \varphi \right)$;\\
  \textbf{A2}\> $\left( \varphi \longrightarrow_{(n)} \left( \psi \longrightarrow_{(n)} \chi \right) \right) \longrightarrow_{(n)} \left( \left( \varphi \longrightarrow_{(n)} \psi \right) \longrightarrow_{(n)} \left( \varphi \longrightarrow_{(n)} \chi \right) \right)$;\\
  \textbf{A3} \> $\left( \neg_{(n)}\psi \longrightarrow_{(n)} \neg_{(n)}\varphi \right) \longrightarrow_{(n)} \left( \left( \neg_{(n)}\psi \longrightarrow_{(n)} \varphi \right) \longrightarrow_{(n)} \psi \right) $; \\
  \textbf{A4} \> $\varphi \longrightarrow_{(n)} \left( \perp_{c_{k}} \longrightarrow_{(n)} \neg_{c_{k}}\varphi \right)$;  \\
  \textbf{A5} \>  $\neg_{c_{k}} \neg_{c_{r}}\varphi \longleftrightarrow_{(n)} \neg_{{c_{k}} \otimes c_{r} } \varphi$;  \\
  \textbf{A6} \> $\neg_{c_{k}}\perp_{c_{r}} \longleftrightarrow_{(n)} \perp_{c_{k} \otimes c_{r}}$; \\
  \textbf{A7} \> $\perp_{c_{k}} \longrightarrow_{(n)} \perp_{c_{r}}$. 
  
\end{tabbing} 

The above axioms hold for every chain $c_{k}$ and $c_{r}$ over $[n]$ with $1 \leq k,r \leq n$. In the case of axiom \textbf{(A7)}, we require that $c_{r}$ be a subchain of $c_{k}$.

\textbf{Rule of inference.}\label{ruleofinference} For all $n\in\mathbb{N}$, the rule of inference for $\textbf{CPN}_{n}$ is \emph{modus ponens}, or simply $MP_{n}$, which states that if one has $\varphi$ and $\varphi \longrightarrow_{(n)} \psi$, then $\psi$ can be deduced.

Since the first three axioms match those of classical propositional calculus, any theorem provable in $\textbf{CPC}$ is also provable in the propositional calculus with multiple negations, using the $n$-connectives defined in definition \ref{def: nconnectives} together with the strong negation $\neg_{(n)}$. The symbol $\vdash_{(n)}$ denotes formal derivation in $\textbf{CPN}_{n}$. Using it, many theorems from classical propositional calculus can be expressed with the generalized connectives of $\textbf{CPN}_{n}$. Our first result is the following proposition.

\begin{proposition}
\label{prop: classicpropositions}
For all $n\in\mathbb{N}$, if $\Sigma$ is a (possibly empty) set of premises of $\textup{\textbf{CPN}}_{n}$, then:

\begin{enumerate}
    \item[\textup{(i)}] $\Sigma, \varphi \vdash_{(n)} \psi$ if and only if $\Sigma \vdash_{(n)} \varphi\longrightarrow_{(n)} \psi$; 
    \item[\textup{(ii)}] $\Sigma, \neg_{(n)} \varphi \vdash_{(n)} \psi$ and $\Sigma, \neg_{(n)}\varphi \vdash_{(n)} \neg_{(n)}\psi $, then $\Sigma \vdash_{(n)} \varphi$; 
    \item[\textup{(iii)}] $\Sigma \vdash_{(n)} \varphi \land_{(n)} \psi$ if and only if $\Sigma \vdash_{(n)} \varphi$ and $\Sigma \vdash_{(n)}  \psi$;  
    \item[\textup{(iv)}] $\Sigma, \neg_{(n)}\varphi \vdash_{(n)} \psi$ if and only if $\Sigma \vdash_{(n)} \varphi \lor_{(n)} \psi$. 
\end{enumerate}

\textup{Item (i) is the Deduction Theorem, and item (ii) is the reductio ad absurdum theorem. In this paper we denote them by $DT_{n}$ and $RT_{n}$, respectively.}

\begin{proof}
The proof is the same as in $\textbf{CPC}$. For further details, see \cite{caicedo}, \cite{Hod}, and \cite{Men}. 
\end{proof}

\begin{corollary}

\label{coro: derivedrulesaxioms}
  For every $n\in\mathbb{N}$ and for all chains $c_{k}$ and $c_{r}$ over $[n]$, the following derived rules hold in $\textup{\textbf{CPN}}_{n}$:  

  \begin{enumerate}
      \item[\textup{(i)}] $\varphi \vdash_{(n)} \perp_{c_{k}} \longrightarrow_{(n)} \neg_{c_{k}}\varphi$ and $\varphi, \perp_{c_{k}} \vdash_{(n)} \neg_{c_{k}}\varphi$; 
      \item[\textup{(ii)}] $\neg_{c_{k}} \neg_{c_{r}}\varphi \vdash_{(n)} \neg_{c_{k} \otimes c_{r}}\varphi$ and $\neg_{c_{k} \otimes c_{r}}\varphi \vdash_{(n)} \neg_{c_{k}} \neg_{c_{r}}\varphi$; 
      \item[\textup{(iii)}] $\neg_{c_{k}}\perp_{c_{r}} \vdash_{(n)} \perp_{c_{k} \otimes c_{r}}$ and $\perp_{c_{k} \otimes c_{r}} \vdash_{(n)} \neg_{c_{k}}\perp_{c_{r}}$; 
      \item[(iv)] $\perp_{c_{k}} \vdash_{(n)} \perp_{c_{r}}$ if $c_{r}$ is a subchain of $c_{k}$.  
  \end{enumerate}

\begin{proof}
Apply $TD_{n}$ and the axioms \textbf{(A4)}, \textbf{(A5)}, \textbf{(A6)}, and \textbf{(A7)}.
\end{proof}

\end{corollary}

\end{proposition}

\begin{proposition}
\label{prop: cknegations}
 For every $n \in \mathbb{N}$ and for all chains $c_{k}$ and $c_{r}$ over $[n]$, the following schemes are derivable in $\textup{\textbf{CPN}}_{n}$: 

\begin{enumerate}
    \item[\textup{(i)}] $\neg_{c_{k}}\neg_{c_{r}} \varphi \longleftrightarrow_{(n)} \neg_{c_{r}}\neg_{c_{k}}\varphi\text{;}$
    \item[\textup{(ii)}] $\neg_{(n)}\neg_{c_{k}} \varphi \longleftrightarrow_{(n)} \neg_{c_{n-k}^{'}}\varphi\text{;}$
    \item[\textup{(iii)}] $\neg_{(n)}\neg_{c_{n-k}^{'}}\varphi \longleftrightarrow_{(n)} \neg_{c_{k}}\varphi\text{;}$
    \item[\textup{(iv)}] $\neg_{c_{k}}\neg_{c_{n-k}^{'}}\varphi \longleftrightarrow_{(n)} \neg_{(n)}\varphi\text{;}$
    \item[\textup{(v)}] $\neg_{c_{k}}\neg_{c_{k}}\varphi \longleftrightarrow_{(n)} \varphi\text{.}$
       
\end{enumerate}

\begin{proof}
Apply axiom \textbf{(A5)} together with Proposition \ref{prop: propietiesofproductofchains}, as well as Remark \ref{rem: conmutativechain} and Notation \ref{not: notationnegationandocntradiction}, to derive the schemes.
\end{proof}
\end{proposition}

\begin{corollary}
\label{coro: derivedrulesnegations}
  For every $n\in\mathbb{N}$ and for all chains $c_{k}$ and $c_{r}$ over $[n]$, the following derived rules hold in $\textup{\textbf{CPN}}_{n}$:   

\begin{enumerate}
    \item[\textup{(i)}] $\neg_{c_{k}} \neg_{c_{r}}\varphi \vdash_{(n)} \neg_{c_{r}}\neg_{c_{k}} \varphi$;
    \item[\textup{(ii)}] $\neg_{(n)} \neg_{c_{k}}\varphi \vdash_{(n)} \neg_{c_{n-k}^{'}}\varphi $ and $\neg_{c_{n-k}^{'}}\varphi \vdash_{(n)} \neg_{(n)}\neg_{c_{k}} \varphi$; 
    \item[\textup{(iii)}] $\neg_{(n)} \neg_{c_{n-k}^{'}}\varphi \vdash_{(n)} \neg_{c_{k}} \varphi$ and $\neg_{c_{k}}\varphi \vdash_{(n)} \neg_{(n)} \neg_{c_{n-k}^{'}}\varphi$; 
    \item[\textup{(iv)}] $\neg_{c_{k}} \neg_{c_{n-k}^{'}}\varphi \vdash_{(n)} \neg_{(n)} \varphi$ and $\neg_{(n)}\varphi \vdash_{(n)} \neg_{c_{k}} \neg_{c_{n-k}^{'}}\varphi$;
    \item[\textup{(v)}] $\neg_{c_{k}} \neg_{c_{k}}\varphi \vdash_{(n)} \varphi$ and $\varphi \vdash_{(n)} \neg_{c_{k}}\neg_{c_{k}}\varphi$.  
\end{enumerate}

\begin{proof}
  Apply $TD_{n}$ to Proposition \ref{prop: cknegations}.  
\end{proof}

\end{corollary}

\begin{proposition}

\label{prop: ckcontradictions}
For every $n\in\mathbb{N}$ and for all chains $c_{k}$ and $c_{r}$ over $[n]$, the following schemes are derivable in $\textup{\textbf{CPN}}_{n}$:

 \begin{enumerate}
     \item[\textup{{(i)}}] $\neg_{c_{r}}\perp_{c_{k}} \longleftrightarrow_{(n)} \neg_{c_{k}}\perp_{c_{r}}$; 
     \item[\textup{(ii)}] $\neg_{(n)}\perp_{c_{k}} \longleftrightarrow_{(n)} \perp_{c_{n-k}^{'}}$;
     \item[\textup{(iii)}] $\neg_{(n)} \perp_{c_{n-k}^{'}} \longleftrightarrow_{(n)} \perp_{c_{k}}$; 
     \item[\textup{(iv)}] $\neg_{c_{k}}\perp_{c_{n-k}^{'}} \longleftrightarrow_{(n)} \perp_{(n)}$; 
     \item[\textup{(v)}] $\neg_{c_{k}}\perp_{c_{k}} \longleftrightarrow_{(n)} \top_{(n)}$. 
 \end{enumerate}

\begin{proof}    
Apply axiom \textbf{(A6)} together with Propositions \ref{prop: chainandcomplementchain} and \ref{prop: propietiesofproductofchains}, as well as Remark \ref{rem: conmutativechain} and Notation \ref{not: notationnegationandocntradiction}.
\end{proof}
 \end{proposition}

\begin{corollary}
\label{coro: derivedrulescontradictions}    
  For every $n\in\mathbb{N}$ and for all chains $c_{k}$ and $c_{r}$ over $[n]$, the following derived rules hold in $\textup{\textbf{CPN}}_{n}$:   
  
 \begin{enumerate}
     \item[\textup{(i)}] $\neg_{c_{r}}\perp_{c_{k}} \vdash_{(n)} \neg_{c_{k}}\perp_{c_{r}}$; 
     \item[\textup{(ii)}] $\neg_{(n)} \perp_{c_{k}} \vdash_{(n)} \perp_{c_{n-k}^{'}}$ and $\perp_{c_{n-k}^{'}} \vdash_{(n)} \neg_{(n)} \perp_{c_{k}}$; 
     \item[\textup{(iii)}] $\neg_{(n)} \perp_{c_{n-k}^{'}} \vdash_{(n)}  \perp_{c_{k}}$ and $\perp_{c_{k}} \vdash_{(n)} \neg_{(n)} \perp_{c_{n-k}^{'}}$; 
     \item[\textup{(iv)}] $\neg_{c_{k}} \perp_{c_{n-k}^{'}} \vdash_{(n)} \perp_{(n)}$ and $\perp_{(n)} \vdash_{(n)} \neg_{c_{k}} \perp_{c_{n-k}^{'}}$; 
     \item[\textup{(v)}] $\neg_{c_{k}} \perp_{c_{k}} \vdash_{(n)} \top_{(n )}$ and $\top_{(n)} \vdash_{(n)} \neg_{c_{k}} \perp_{c_{k}}$.  
 \end{enumerate} 

 \begin{proof}
   Apply $TD_{n}$ to Proposition \ref{prop: ckcontradictions}.    
 \end{proof}
\end{corollary}

\begin{proposition}
\label{prop: factorizationprop}
For every $n\in\mathbb{N}$, if $c_{k}$ and $c_{r}$ are chains over $[n]$, then there exists a chain $c_{s}$ over $[n]$ such that the following schemes are derivable in $\textup{\textbf{CPN}}_{n}$:

\begin{enumerate}
    \item[\textup{(i)}] $\neg_{c_{k}}\varphi \longleftrightarrow_{(n)} \neg_{c_{s}}\neg_{c_{r}} \varphi$;
    \item[\textup{(ii)}] $\perp_{c_{k}} \longleftrightarrow_{(n)} \neg_{c_{s}}\perp_{c_{r}}$.  
\end{enumerate}
\begin{proof}
 For both items, it suffices to take $c_{s}=c_{k}\cdot c_{r}$, apply axioms \textbf{(A5)} and \textbf{(A6)}, and then use Proposition \ref{prop: propietiesofproductofchains}.
\end{proof}
\end{proposition}

Since $\textup{\textbf{CPN}}_{n}$ behaves the same as $\textup{\textbf{CPC}}$ thanks to axioms \textbf{(A1)}, \textbf{(A2)}, and \textbf{(A3)}, any formula built from the connectives $\longrightarrow_{(n)}$, $\land_{(n)}$, $\lor_{(n)}$, and $\neg_{(n)}$ needs no proof; it is the same as in classical propositional calculus. Therefore, we will use these theorems freely. From now on, some proof lines will invoke these theorems; we will mark this with $\textup{CF}_{n}$ (classic formula in $\textup{\textbf{CPN}}_{n}$ ), and the specific theorem used will be listed in the notes at the end of the paper.

\begin{proposition}
\label{prop: ckfirstpropositions}
For all $n\in\mathbb{N}$ and for every chain $c_{k}$ over $[n]$, the following schemes are derivable in $\textup{\textbf{CPN}}_{n}$:

\begin{enumerate}
\item[\textup{(i)}] $\varphi \land_{(n)} \neg_{c_{k}} \varphi \longrightarrow_{(n)} \perp_{c_{k}}$;
\item[\textup{(ii)}] $\perp_{c_{n-k}^{'}} \longrightarrow_{(n)} \varphi \lor_{(n)} \neg_{c_{k}}\varphi$;
\item[\textup{(iii)}] $\neg_{c_{k}}\varphi \longrightarrow_{(n)} \left( \perp_{c_{k}} \longrightarrow_{(n)} \varphi \right)$;
\item[\textup{(iv)}] $\neg_{c_{k}}\varphi \longrightarrow_{(n)} \left( \perp_{c_{n-k}^{'}} \longrightarrow_{(n)} \neg_{(n)}\varphi \right)$;
\item[\textup{(v)}] $\neg_{c_{k}}\varphi \longrightarrow_{(n)} \left( \varphi \longleftrightarrow_{(n)} \perp_{c_{k}} \right)$; 
\item[\textup{(vi)}] $\neg_{c_{k}}\varphi \land_{(n)} \neg_{c_{n-k}^{'}}\varphi \longrightarrow_{(n)} \perp_{(n)}$; 
\item[\textup{(vii)}] $\perp_{c_{k}} \land_{(n)} \perp_{c_{n-k}^{'}} \longrightarrow_{(n)} \perp_{(n)}\text{.}$
\end{enumerate}

\begin{proof} 

\item[(i)] Let us show that $\varphi \land_{(n)} \neg_{c_{k}} \varphi \vdash_{(n)} \perp_{c_{k}}$

 {\footnotesize{

 \begin{longlogicproof}{1}
     \varphi \land_{(n)} \neg_{c_{k}}\varphi & Premise \\
     
    \varphi  & $\mathrm{CF}_{n}$\footnote{Use the theorem $\varphi \land_{(n)} \psi \vdash_{(n)} \varphi$.}  to  (1)  \\
    
     \neg_{c_{k}} \varphi & $\mathrm{CF}_{n}$\footnote{Use the theorem $\varphi \land_{(n)} \psi \vdash_{(n)} \psi$.} to (1) \\ 
     
     \perp_{c_{n-k}^{'}} \longrightarrow_{(n)} \neg_{c_{n-k}^{'}} \varphi  & Corollary \ref{coro: derivedrulesaxioms}-(i) to (2)    \\
         
     \neg_{c_{k}}\varphi  \longrightarrow_{(n)} \perp_{c_{k}} & $\mathrm{CF}_{n}$\footnote{Use the theorem $\psi \longrightarrow_{(n)} \varphi \vdash_{(n)} \neg_{(n)}\varphi \longrightarrow_{(n)} \neg_{(n)}\psi $ .} + Corollary \ref{coro: derivedrulesnegations}, \ref{coro: derivedrulescontradictions} to (4) \\

     \perp_{c_{k}} & $MP(5,3)$
         
\end{longlogicproof}

}}
Using $TD_{n}$, we obtain the desired result.
\vspace{0.5cm}
\item[(ii)] Let us show that $\perp_{c_{n-k}^{'}}, \neg_{(n)}\varphi \vdash_{(n)} \neg_{c_{k}}\varphi$
{\footnotesize{

 \begin{longlogicproof}{1}
     \perp_{c_{n-k}^{'}} & Premise \\
     
     \neg_{(n)}\varphi  & Premise  \\
    
    \perp_{c_{n-k}^{'}} \longrightarrow_{(n)} \neg_{c_{n-k}^{'}}\neg_{(n)}\varphi & Corollary \ref{coro: derivedrulesaxioms}-(i) to (2) \\ 
     
     \neg_{c_{n-k}^{'}}\neg_{(n)}\varphi & $MP_{n}(3,2)$    \\
         
     \neg_{c_{k}}\varphi & Corollary \ref{coro: derivedrulesnegations}-(iii)  to (4) 
         
\end{longlogicproof}
}}
Using Proposition \ref{prop: classicpropositions}-(iv), we obtain the desired result.
\vspace{0.5cm}
\item[(iii)] Let us show that $\neg_{c_{k}}\varphi, \perp_{c_{k}} \vdash_{(n)} \varphi$

 {\footnotesize{

 \begin{longlogicproof}{1}
     \neg_{c_{k}}\varphi & Premise \\
     
     \perp_{c_{k}}  & Premise  \\
    
    \perp_{c_{k}} \longrightarrow_{(n)} \neg_{c_{k}}\neg_{c_{k}}\varphi & Corollary \ref{coro: derivedrulesaxioms}-(i) to (1) \\ 
     
     \neg_{c_{k}}\neg_{c_{k}}\varphi & $MP_{n}(3,2)$    \\
         
     \varphi & Corollary \ref{coro: derivedrulesnegations}-(v) to (4) 
         
\end{longlogicproof}
}}
Using $TD_{n}$, we obtain the desired result.
\vspace{0.5cm}
\item[(iv)] Let us show that $\neg_{c_{k}}\varphi, \perp_{c_{n-k}^{'}} \vdash_{(n)}\neg_{(n)} \varphi$

{\footnotesize{

 \begin{longlogicproof}{1}
     \neg_{c_{k}}\varphi & Premise \\
     
     \perp_{c_{n-k}^{'}}  & Premise  \\
    
    \perp_{c_{n-k}^{'}} \longrightarrow_{(n)} \neg_{c_{n-k}^{'}}\neg_{c_{k}}\varphi & Corollary \ref{coro: derivedrulesaxioms}-(i) to (1) \\ 
     
     \neg_{c_{n-k}^{'}}\neg_{c_{k}}\varphi & $MP_{n}(3,2)$    \\
         
     \neg_{(n)}\varphi & Corollary \ref{coro: derivedrulesnegations}-(iv) to (4) 
         
\end{longlogicproof}
}}
Using $TD_{n}$, we obtain the desired result.
\vspace{0.5cm}
\item[(v)] Let us show that $ \neg_{c_{k}}\varphi \vdash_{(n)} \varphi \longleftrightarrow_{(n)} \perp_{c_{k}}$

{\footnotesize{

 \begin{longlogicproof}{1}
     \neg_{c_{k}}\varphi & Premise \\
     
      \perp_{c_{k}} \longrightarrow_{(n)} \varphi & Item (iii) to (1)  \\
    
    \perp_{c_{n-k}^{'}} \longrightarrow_{(n)} \neg_{(n)}\varphi & Item (iv) to (1) \\ 
     
    \varphi \longrightarrow_{(n)} \perp_{c_{k}} & $\text{CF}_{n}$\footnote{Use the theorem $\neg_{(n)}\psi \longrightarrow_{(n)} \neg_{(n)}\varphi \vdash_{(n)} \varphi \longrightarrow_{(n)} \psi $.} + Corollary \ref{coro: derivedrulesnegations}-(v), \ref{coro: derivedrulescontradictions}-(iii) to (3)   \\
         
    \varphi \longleftrightarrow_{(n)} \perp_{c_{k}} & $\text{CF}_{n}$\footnote{Use the theorem $\varphi \longrightarrow_{(n)} \psi, \psi \longrightarrow_{(n)} \varphi \vdash_{(n)} \varphi \longleftrightarrow_{(n)}\psi$.} to (2,4)        
\end{longlogicproof}
}}

Using $TD_{n}$, we obtain the desired result.
\vspace{0.5cm}
\item[(vi)] Let us show that $ \neg_{c_{k}}\varphi \land_{(n)} \neg_{c_{n-k}^{'}}\varphi \vdash_{(n)} \perp_{(n)}$

{\footnotesize{

 \begin{longlogicproof}{1}
     \neg_{c_{k}}\varphi \land_{(n)} \neg_{c_{n-k}^{'}}\varphi & Premise \\
     
     \neg_{c_{k}}\varphi  & $\textup{CF}_{n}$\footnote{See note 1. } to (1)  
 \\
    
    \neg_{c_{n-k}^{'}}\varphi & $\textup{CF}_{n}$\footnote{See note 2.} to (1) \\ 
     
     \neg_{(n)}\neg_{c_{k}}\varphi & Corollary \ref{coro: derivedrulesnegations}-(ii) to (3)   \\
         
    \perp_{(n)} & $\textup{CF}_{n}$\footnote{Use the theorem $\varphi,\neg_{(n)}\varphi \vdash_{(n)} \perp_{(n)}\text{.}$} to (2,4)
         
\end{longlogicproof}
}}
Using $TD_{n}$, we obtain the desired result.
\vspace{0.5cm}

\item[(vii)] The proof is identical to that of item (v), except that here we use Corollary \ref{coro: derivedrulescontradictions}.

\end{proof}
\end{proposition}

\begin{corollary}
\label{coro: derivedrulesfirstpropositions}
For every $n\in\mathbb{N}$ and for every chain $c_{k}$ over $[n]$, the following derived rules hold in $\textup{\textbf{CPN}}_{n}$: 

\begin{enumerate}
    \item[\textup{(i)}] $\varphi \land_{(n)} \neg_{c_{k}}\varphi \vdash_{(n)} \perp_{c_{k}}$ and $\varphi, \neg_{c_{k}}\varphi \vdash_{(n)} \perp_{c_{k}}$; 
    \item[\textup{(ii)}] $\perp_{c_{n-k}^{'}} \vdash_{(n)} \varphi \lor_{(n)} \neg_{c_{k}}\varphi$ and $\perp_{c_{n-k}^{'}}, \neg_{(n)} \varphi \vdash_{(n)} \neg_{c_{k}}\varphi$;  
    \item[\textup{(iii)}] $\neg_{c_{k}}\varphi \vdash_{(n)} \perp_{c_{k}} \longrightarrow_{(n)} \varphi$ and $\neg_{c_{k}}\varphi, \perp_{c_{k}} \vdash_{(n)} \varphi$; 
    \item[\textup{(iv)}] $\neg_{c_{k}}\varphi \vdash_{(n)} \perp_{c_{n-k}^{'}} \longrightarrow_{(n)} \neg_{(n)} \varphi$ and $\neg_{c_{k}}\varphi, \perp_{c_{n-k}^{'}} \vdash_{(n)} \neg_{(n)} \varphi$; 
    \item[\textup{(v)}] $\neg_{c_{k}}\varphi \land_{(n)} \neg_{c_{n-k}^{'}}\varphi \vdash_{(n)} \perp_{(n)}$ and $\neg_{c_{k}}\varphi, \neg_{c_{n-k}^{'}}\varphi \vdash_{(n)} \perp_{(n)}$;  
    \item[\textup{(vi)}] $\perp_{c_{k}} \land_{(n)} \perp_{c_{n-k}^{'}} \vdash_{(n)} \perp_{(n)}$ and $\perp_{c_{k}}, \perp_{c_{n-k}^{'}} \vdash_{(n)} \perp_{(n)}$.   
\end{enumerate}

\begin{proof}
  Use Proposition \ref{prop: classicpropositions} and Proposition \ref{prop: ckfirstpropositions}.  
\end{proof}
 
\end{corollary}

\begin{proposition}
\label{prop: cksecondpropositions}
 For all $n\in\mathbb{N}$ and for every chain $c_{k}$ over $[n]$, the following schemes are derivable in $\textup{\textbf{CPN}}_{n}$:

\begin{enumerate}
    \item[\textup{(i)}] $\perp_{c_{k}}\longrightarrow_{(n)} \neg_{c_{n-k}^{'}}\left( \varphi \land_{(n)} \neg_{c_{n-k}^{'}}\varphi \right)$;
    \item[\textup{(ii)}] $\neg_{c_{k}}\varphi \land_{(n)} \neg_{c_{k}}\psi \longrightarrow_{(n)} \neg_{c_{k}}\left( \varphi \land_{(n)} \psi \right)$;
    \item[\textup{(iii)}] $\neg_{c_{k}}\left(\varphi \lor_{(n)} \psi \right)\longrightarrow_{(n)} \neg_{c_{k}}\varphi \lor_{(n)} \neg_{c_{k}}\psi$;
    \item[\textup{(iv)}] $\neg_{c_{k}}\left( \varphi \land_{(n)} \psi \right) \longrightarrow_{(n)} \neg_{c_{k}}\varphi \lor_{(n)} \neg_{c_{k}} \psi$;
    \item[\textup{(v)}] $\neg_{c_{k}}\varphi \land_{(n)} \neg_{c_{k}}\psi \longrightarrow_{(n)} \neg_{c_{k}} \left( \varphi \lor_{(n)} \psi \right)$.

\end{enumerate}

\begin{proof}

    \item[(i)] Let us show that $\perp_{c_{k}}, \ \neg_{c_{k}}\left( \varphi \land_{(n)} \neg_{c_{n-k}^{'}}\varphi \right) \ \vdash_{(n)} \ \perp_{(n)}$ \footnote{One equivalent formulation of the reductio ad absurdum theorem states that if $\Sigma, \neg_{(n)}\varphi \vdash_{(n)} \perp_{(n)}$, then $\Sigma \vdash_{(n)} \varphi$.}

 {\footnotesize{

 \begin{longlogicproof}{1}
     \perp_{c_{k}}& Premise \\
     \neg_{c{_k}}\left( \varphi \land_{(n)} \neg_{c_{n-k}^{'}}\varphi \right) & Premise  \\
    \perp_{c_{k}} \longrightarrow_{(n)} \varphi \land_{(n)} \neg_{c_{n-k}^{'}}\varphi & Corollary \ref{coro: derivedrulesfirstpropositions}-(iii) to (2) \\ 
         \varphi \land_{(n)} \neg_{c_{n-k}^{'}}\varphi & $MP_{n}(3,1)$    \\      
     \perp_{c_{n-k}^{'}} & Corollary \ref{coro: derivedrulesfirstpropositions}-(i) to (4)  \\
     \perp_{(n)} & Corollary \ref{coro: derivedrulesfirstpropositions}-(vi) to (1,5)     
\end{longlogicproof}
}}
By the reductio ad absurdum theorem, we have $\perp_{c_{k}} \vdash_{(n)} \neg_{c_{n-k}^{'}} \left( \varphi \land_{(n)} \neg_{c_{n-k}^{'}}\varphi\right)$; applying $TD_{n}$ gives the desired result.
\vspace{0.5cm}
\item[(ii)] Let us show that $\{ \neg_{c_{k}}\varphi, \neg_{c_{k}}\psi\}, \neg_{c_{n-k}^{'}}\left( \varphi \land_{(n)} \psi \right) \vdash_{(n)} \perp_{(n)}$ 

{\footnotesize{

 \begin{longlogicproof}{1}
     \neg_{c_{k}}\varphi & Premise \\
     
     \neg_{c_{k}}\psi & Premise  \\
    
    \neg_{c_{n-k}^{'}}\left( \varphi \land_{(n)} \psi \right) & Premise \\ 
     
    \perp_{c_{k}} \longrightarrow_{(n)} \varphi & Corollary \ref{coro: derivedrulesfirstpropositions}-(iii) to (1)  \\
         
     \perp_{c_{k}} \longrightarrow_{(n)} \psi & Corollary \ref{coro: derivedrulesfirstpropositions}-(iii) to (2) \\

     \perp_{c_{k}} \longrightarrow_{(n)} \varphi \land_{(n)} \psi & $\text{CF}_{n}$\footnote{Use the theorem $\varphi \longrightarrow_{(n)} \psi, \varphi\longrightarrow_{(n)} \chi \vdash_{(n)} \varphi \longrightarrow_{(n)} \psi \land_{(n)} \chi \text{.}$} to (4,5) \\

     \perp_{c_{k}} \longrightarrow_{(n)} \neg_{(n)}\left( \varphi \land_{(n)} \psi \right)  & Corollary \ref{coro: derivedrulesfirstpropositions}-(iv) to (3) \\

     \perp_{c_{n-k}^{'}} & $\text{CF}_{n}$\footnote{Use the theorem $\varphi\longrightarrow_{(n)}\psi, \varphi\longrightarrow_{(n)}\neg_{(n)}\psi \vdash_{(n)} \neg_{(n)}\varphi\text{.}$} to (6,7) + Corollary \ref{coro: derivedrulescontradictions}-(ii) \\

     \perp_{c_{n-k}^{'}} \longrightarrow_{(n)} \neg_{(n)}\varphi & Corollary \ref{coro: derivedrulesfirstpropositions}-(iv) to (1) \\

     \neg_{(n)}\varphi & $MP_{n}(9,8)$ \\

     \perp_{c_{n-k}^{'}} \longrightarrow_{(n)} \varphi \land_{(n)} \psi & Corollary \ref{coro: derivedrulesfirstpropositions}-(iii) to (3) \\

     \varphi \land_{(n)} \psi & $MP_{n}(11,8)$ \\

     \varphi & $\text{CF}_{n}$\footnote{See note 1.} to (12) \\

     \perp_{(n)} & $\text{CF}_{n}$\footnote{See note 9.} to (10,13)        
\end{longlogicproof}
}}

By the reductio ad absurdum theorem, we have $\neg_{c_{k}}\varphi, \neg_{c_{k}}\psi \vdash_{(n)} \neg_{c_{k}}\left( \varphi \land_{(n)} \psi \right)$; applying item (iii) of Proposition \ref{prop: classicpropositions} and $TD_{n}$ gives the desired result.
\vspace{0.5cm}

\item[(iii)]  Let us show that $ \{ \neg_{c_{k}}\left( \varphi \lor_{(n)} \psi \right), \neg_{c_{n-k}^{'}}\varphi \}, \neg_{c_{n-k}^{'}}\psi \vdash_{(n)} \perp_{(n)} $

 {\footnotesize{

 \begin{longlogicproof}{1}
     \neg_{c_{k}} \left( \varphi \lor_{(n)} \psi \right) & Premise \\
     
    \neg_{c_{n-k}^{'}}\varphi & Premise  \\
    
    \neg_{c_{n-k}^{'}}\psi & Premise \\ 
     
     \perp_{c_{k}} \longrightarrow_{(n)} \neg_{(n)}\varphi & Corollary \ref{coro: derivedrulesfirstpropositions}-(iv) to (2)   \\
         
     \perp_{c_{k}} \longrightarrow_{(n)} \varphi \lor_{(n)} \psi & Corollary \ref{coro: derivedrulesfirstpropositions}-(iii) to (1) \\

     \perp_{c_{k}} \longrightarrow_{(n)} \psi &  $\text{CF}_{n}$\footnote{Use the theorem $\varphi \longrightarrow_{(n)} \psi \lor \chi, \varphi\longrightarrow_{(n)} \neg_{(n)}\psi \vdash_{(n)} \varphi \longrightarrow_{(n)} \chi \text{.}$} to (5,4) \\

     \perp_{c_{k}} \longrightarrow_{(n)}  \neg_{(n)} \psi & Corollary \ref{coro: derivedrulesfirstpropositions}-(iv) to (3) \\

     \perp_{c_{n-k}^{'}} & $\text{CF}_{n}$\footnote{See note 9.} to (6,7) + Corollary \ref{coro: derivedrulescontradictions}-(ii) \\

     \perp_{c_{n-k}^{'}} \longrightarrow_{(n)} \neg_{(n)} \left( \varphi \lor_{(n)} \psi \right) & Corollary \ref{coro: derivedrulesfirstpropositions}-(iv) to (1) \\

     \neg_{(n)}\left( \varphi \lor_{(n)} \psi \right) & $MP_{n}(9,8)$ \\

     \neg_{(n)}\varphi & $\text{CF}_{n}$\footnote{Use the theorem $\neg_{(n)}\left( \varphi \lor_{(n)} \psi \right) \vdash_{(n)} \neg_{(n)}\varphi$. } to (10) \\ 

     \perp_{c_{n-k}^{'}} \longrightarrow_{(n)} \varphi & Corollary \ref{coro: derivedrulesfirstpropositions}-(iii) to (2) \\

     \varphi & $MP_{n}(12,8)$ \\

     \perp_{(n)} & $\text{CF}_{n}$\footnote{See note 8.} to (13,11)
         
\end{longlogicproof}
}}
By the reductio ad absurdum theorem we have $\neg_{c_{k}}\left( \varphi \lor_{(n)} \psi \right), \neg_{c_{n-k}^{'}}\varphi \vdash_{(n)} \neg_{c_{k}} \psi$.
Applying item (iv) of Proposition \ref{prop: classicpropositions} we have $\neg_{c_{k}} \left( \varphi \lor_{(n)} \psi \right) \vdash_{(n)} \neg_{c_{k}}\varphi \lor_{(n)} \neg_{c_{k}} \psi$ and using $TD_{n}$ we obtain the desired result.
\vspace{0.5cm}

\item[(iv)]  Let us show that $\{ \neg_{c_{k}}\left( \varphi \land_{(n)} \psi \right), \neg_{c_{n-k}^{'}}\varphi\}, \neg_{c_{n-k}^{'}}\psi \vdash_{(n)} \perp_{(n)} $

{\footnotesize{

 \begin{longlogicproof}{1}
     \neg_{c_{k}}\left( \varphi \land_{(n)} \psi \right) & Premise \\
     
    \neg_{c_{n-k}}^{'}\varphi & Premise  \\
    
   \neg_{c_{n-k}^{'}}\psi &  Premise \\ 
     
     \neg_{c_{n-k}^{'}}\varphi \land_{(n)} \neg_{c_{n-k}^{'}}\psi & $\text{CF}_{n}$\footnote{Use the theorem $\varphi, \psi \vdash_{(n)} \varphi \land_{(n)} \psi$.} to (2,3)  \\
         
     \neg_{c_{n-k}^{'}} \left( \varphi \land_{(n)} \psi \right) & Proposition \ref{prop: cksecondpropositions} + $TD_{n}$ to (4) \\

     \perp_{(n)} & Corollary \ref{coro: derivedrulesfirstpropositions}-(v) to (1,5)
         
\end{longlogicproof}
}}

By the reductio ad absurdum theorem we have $\neg_{c_{k}}\left( \varphi \land_{(n)} \psi \right), \neg_{c_{n-k}^{'}}\varphi \vdash_{(n)} \neg_{c_{k}} \psi$.
Applying item (iv) of Proposition \ref{prop: classicpropositions} we have $\neg_{c_{k}} \left( \varphi \land_{(n)} \psi \right) \vdash_{(n)} \neg_{c_{k}}\varphi \lor_{(n)} \neg_{c_{k}} \psi$ and using $TD_{n}$ we obtain the desired result.
\vspace{0.5cm}
\item[(v)]  Let us show that $ \neg_{c_{k}}\varphi \land_{(n)} \neg_{c_{k}}\psi, \neg_{c_{n-k}^{'}}\left( \varphi \lor_{(n)} \psi \right) \vdash_{(n)} \perp_{(n)}$

{\footnotesize{

 \begin{longlogicproof}{1}
     \neg_{c_{k}}\varphi \land_{(n)} \neg_{c_{k}}\psi & Premise \\
     
    \neg_{c_{n-k}^{'}}\left( \varphi \lor_{(n)} \psi \right) & Premise  \\
    
   \neg_{c_{k}}\varphi &  $\text{CF}_{n}$\footnote{See note 1.} to (1) \\ 
     
     \neg_{c_{k}}\psi & $\text{CF}_{n}$\footnote{See note 2.} to (1)  \\
         
     \perp_{c_{k}} \longrightarrow_{(n)} \neg_{(n)} \left( \varphi \lor_{(n)} \psi \right) & Corollary \ref{coro: derivedrulesfirstpropositions}-(iv) to (2) \\

     \perp_{c_{k}}\longrightarrow_{(n)} \neg_{(n)}\varphi \land_{(n)} \neg_{(n)}\psi & $\text{CF}_{n}$\footnote{Use the theorem $\varphi \longrightarrow_{(n)} \neg_{(n)} \left( \psi \lor_{(n)} \chi \right) \vdash_{(n)} \varphi \longrightarrow_{(n)} \neg_{(n)}\psi \land_{(n)} \neg_{(n)} \chi\text{.}$} to (5) \\

     \left( \perp_{c_{k}}  \longrightarrow_{(n)} \neg_{(n)}\varphi \right) \land_{(n)} \left( \perp_{c_{k}} \longrightarrow_{(n)} \neg_{(n)}\psi \right)  & $\text{CF}_{n}$\footnote{Use the theorem $\varphi \longrightarrow_{(n)} \psi \land_{(n)} \chi \vdash_{(n)} \left( \varphi \longrightarrow_{(n)} \psi \right) \land_{(n)} \left( \varphi \longrightarrow_{(n)} \chi  \right) \text{.}$} to (6) \\

     \perp_{c_{k}} \longrightarrow_{(n)} \neg_{(n)}\varphi & $\text{CF}_{n}$\footnote{See note 1.} to (7) \\

    \perp_{c_{k}} \longrightarrow_{(n)} \neg_{(n)}\psi & $\text{CF}_{n}$\footnote{See note 2.} to (7) \\

    \perp_{c_{k}} \longrightarrow_{(n)} \varphi & Corollary \ref{coro: derivedrulesfirstpropositions}-(iii) to (3) \\

    \perp_{c_{n-k}^{'}} & $\text{CF}_{n}$\footnote{See note 12.} to (8,10) + Corollary \ref{coro: derivedrulescontradictions}-(ii) \\

    \perp_{c_{n-k}^{'}} \longrightarrow_{(n)} \neg_{(n)}\varphi & Corollary \ref{coro: derivedrulesfirstpropositions}-(iv) to (3) \\

    \neg_{(n)}\varphi & $MP_{n}(12,11)$ \\

    \perp_{c_{n-k}^{'}} \longrightarrow_{(n)} \neg_{(n)}\psi & Corollary \ref{coro: derivedrulesfirstpropositions}-(iv) to (4)\\

    \neg_{(n)}\psi & $MP_{n}(14,11)$ \\

    \perp_{c_{n-k}^{'}} \longrightarrow_{(n)} \varphi \lor_{(n)} \psi & Corollary \ref{coro: derivedrulesfirstpropositions}-(iii) to (2) \\

    \varphi \lor_{(n)} \psi & $MP_{n}(16,11)$ \\

    \varphi & $\text{CF}_{n}$\footnote{Use the theorem $\varphi \lor_{(n)} \psi, \neg_{(n)}\psi \vdash_{(n)} \varphi \text{.}$} to (17,15) \\

    \perp_{(n)} & $\text{CF}_{n}$\footnote{See note 9. } to (13,18)
  
\end{longlogicproof}
}}

By the reductio ad absurdum theorem we have $\neg_{c_{k}}\varphi \land_{(n)} \neg_{c_{k}}\psi \vdash_{(n)} \neg_{c_{k}} \left( \varphi \lor_{(n)} \psi \right)$. Applying $TD_{n}$ we obtain the desired result.
\end{proof}

\end{proposition}

\begin{theorem}
\label{teo: fundamental}
For every $n\in\mathbb{N}$, if $c_{k}$ and $c_{s}$ are chains over $[n]$ and $c_{r}$ is the chain of symbols common to $c_{k}$ and $c_{s}$, then:

$$ \vdash_{(n)} \neg_{c_{k}}\varphi \longrightarrow_{(n)} \left( \neg_{c_{s}}\psi \longrightarrow_{(n)} \neg_{c_{s}\otimes c_{r}}\left( \varphi \longrightarrow_{(n)} \psi \right) \right)\text{.} $$

\begin{proof}

Let $c_{p}=c_{s}\otimes c_{r}$, we show that $\neg_{c_{k}}\varphi, \neg_{c_{s}}\psi \vdash_{(n)} \neg_{c_{p}} \left( \varphi \longrightarrow_{(n)} \psi \right)$, to this end we first show that $\neg_{c_{k}}\varphi, \neg_{c_{s}}\psi, \neg_{c_{n-p}^{'}}\left( \varphi \longrightarrow_{(n)} \psi \right) \vdash_{(n)} \perp_{(n)}$, and applying $RA_{n}$ gives the desired result.
{\footnotesize{

 \begin{longlogicproof}{1}
      \neg_{c_{k}}\varphi & Premise \\  
    \neg_{c_{s}}\psi & Premise  \\
    \neg_{c_{n-p}^{'}}\left( \varphi \longrightarrow_{(n)} \psi \right) & Premise \\ 
     
    \perp_{c_{p}} \longrightarrow_{(n)} \neg_{(n)}\left( \varphi \longrightarrow_{(n)} \psi \right) & Corollary \ref{coro: derivedrulesfirstpropositions}-(iv) to (3)  \\
         
    \perp_{c_{s}} \longrightarrow_{(n)} \perp_{c_{p}} & (A7), since $c_{p}$ is a subchain of $c_{s}$ \\

     \perp_{c_{s}} \longrightarrow_{(n)} \neg_{(n)} \left( \varphi \longrightarrow_{(n)}\psi \right) & $\text{CF}_{n}$\footnote{Use the theorem $\varphi \longrightarrow_{(n)} \psi, \psi \longrightarrow_{(n)} \chi \vdash_{(n)} \varphi \longrightarrow_{(n)} \chi$.} to (5,4) \\

      \perp_{c_{s}} \longrightarrow_{(n)} \psi  & Corollary \ref{coro: derivedrulesfirstpropositions}-(iii) to (2) \\

      \perp_{c_{s}} \longrightarrow_{(n)} \left( \varphi \longrightarrow_{(n)} \psi \right) & $\text{CF}_{n}$\footnote{Use the theorem $\varphi\longrightarrow_{(n)}\psi \vdash_{(n)} \varphi\longrightarrow_{(n)} \left( \chi \longrightarrow_{(n)} \psi \right)$.} to (7) \\

      \perp_{c_{n-s}^{'}} & $\text{CF}_{n}$\footnote{See note 12.} to (6,8) + Corollary \ref{coro: derivedrulescontradictions}-(ii) \\

      \perp_{c_{n-s}^{'}} \longrightarrow_{(n)} \neg_{(n)}\psi & Corollary \ref{coro: derivedrulesfirstpropositions}-(iv) to (2) \\

      \neg_{(n)}\psi & $MP_{n}(10,9)$ \\

      \perp_{c_{n-p}^{'}} \longrightarrow_{(n)} \perp_{c_{k}}   & (A7), since $c_{k}$ is a subchain of $c_{n-p}^{'}$ \\
      \perp_{c_{k}} \longrightarrow_{(n)} \varphi & Corollary \ref{coro: derivedrulesfirstpropositions}-(iii) to (1) \\

      \perp_{c_{n-p}^{'}} \longrightarrow_{(n)} \varphi & $\text{CF}_{n}$\footnote{See note 29.} to (12,13) \\

       \perp_{c_{n-p}^{'}} \longrightarrow_{(n)} \left( \varphi \longrightarrow_{(n)} \psi  \right)  & Corollary \ref{coro: derivedrulesfirstpropositions}-(iii) to (3) \\

       \perp_{c_{n-p}^{'}} \longrightarrow_{(n)} \psi & $\text{CF}_{n}$\footnote{Use the theorem $\varphi\longrightarrow_{(n)} \left( \psi \longrightarrow_{(n)} \chi \right), \varphi\longrightarrow_{(n)} \psi \vdash_{(n)} \varphi\longrightarrow_{(n)} \chi$.} to (15,14) \\

       \perp_{c_{p}} & $\text{CF}_{n}$\footnote{Use the theorem $\varphi \longrightarrow_{(n)} \psi, \neg_{(n)}\psi \vdash_{(n)} \neg_{(n)}\varphi$.} to (16,11) + Corollary \ref{coro: derivedrulescontradictions}-(iii)  \\

       \neg_{(n)} \left( \varphi \longrightarrow_{(n)} \psi \right) & $MP_{n}(17,4)$ \\

       \varphi  & $\text{CF}_{n}$\footnote{Use the theorem $\neg_{(n)} \left( \varphi \longrightarrow_{(n)} \psi \right) \vdash_{(n)} \varphi\text{.}$} to (18) \\

       \perp_{c_{k}} & Corollary \ref{coro: derivedrulesfirstpropositions}-(i) to (19,1) \\

       \perp_{c_{p}} \longrightarrow_{(n)} \neg_{c_{p}}\perp_{c_{k}} & Corollary \ref{coro: derivedrulesaxioms}-(i) to (20) \\

       \neg_{c_{p}}\perp_{c_{k}} & $MP_{n}(21,17)$ \\

       \perp_{c_{p}\otimes c_{k}} & Corollary \ref{coro: derivedrulesaxioms}-(iii) to (22) \\

       \perp_{c_{p} \otimes c_{k}} \longrightarrow_{(n)} \perp_{c_{s}} &(A7), since $c_{s}$ is subchain of $c_{p}\otimes c_{k}$ \\

       \perp_{c_{s}} & $MP_{n}(24,23)$ \\

       \perp_{(n)} & Corollary \ref{coro: derivedrulesfirstpropositions}-(vi) to  (25,9)

\end{longlogicproof}
}}

Using $TD_{n}$, the proof of the theorem is concluded.
\end{proof}
\end{theorem}

\begin{corollary}
 For every $n\in\mathbb{N}$, if $c_{k}$ and $c_{s}$ are chains over $[n]$, then:

\begin{enumerate}
\item[\textup{(i)}] If $c_{s}$ is a subchain of $c_{k}$, then $\vdash_{(n)} \neg_{c_{k}}\varphi \longrightarrow_{(n)} \left( \neg_{c_{s}}\psi \longrightarrow \left( \varphi \longrightarrow \psi \right) \right)$;

\item[\textup{(ii)}] If $c_{k}$ is a subchain of $c_{s}$, then $\vdash_{(n)} \neg_{c_{k}}\varphi \longrightarrow_{(n)} \left( \neg_{c_{s}}\psi \longrightarrow \neg_{c_{k} \otimes c_{s}} \left( \varphi \longrightarrow \psi \right) \right)$;

\item[\textup{(iii)}] If $c_{k}$ and $c_{s}$ have no symbols in common, then $\vdash_{(n)} \neg_{c_{k}}\varphi \longrightarrow_{(n)} \left( \neg_{c_{s}}\psi \longrightarrow \neg_{c_{s}} \left( \varphi \longrightarrow \psi \right) \right)$.
\end{enumerate}   

\begin{proof}
Each result follows by an appropriate application of Theorem \ref{teo: fundamental} regarding the conditions on the chains.
\end{proof}

\end{corollary}

\section{Semantics for $\textbf{CPN}_{n}$}
In this section we develop the semantics of multiple worlds for $\textbf{CPN}_{n}$. Since we work with several negations, each one corresponds to a distinct world where its action on the propositions it negates is validated, while in other worlds it is ignored. Negations interact only within their own world. By a world we mean a set $\mathcal{M}_{i} = \{ T_{i}, F_{i} \}$, where $T_{i}$ stands for truth and $F_{i}$ for falsity, and where the subscript $i$ ranges from $1$ to $n$.

\begin{definition}
\label{def: ivaluation}
Let $\{ \mathcal{M}_{i} \}_{1 \leq i \leq n}$ be a collection of $n$ worlds. For each $1 \leq i \leq n$, we define the $i$-valuation as the function $v_{i} \colon \Omega_{n} \to \mathcal{M}_{i}\text{.}$
\end{definition}

\begin{definition}
 Each $i$-valuation extends to $\mathcal{F}\left( \Omega_{n} \right)$ as follows. For each chain $c_{k}$ over $[n]$, define $\overline{v_{i}} \colon \mathcal{F}\left( \Omega_{n} \right) \to \mathcal{M}_{i}$ by:
\begin{enumerate}
    \item[1.] $\overline{v_{i}} \left( \varphi \right)= v_{i}(\varphi)$ if $\varphi\in P_{n}$.
    \item[2.] $ \overline{v_{i}}\left( \perp_{c_{k}} \right)=\begin{cases}
    F_{i} \ , & \text{if $i$ is a symbol of $c_{k}$}\\
    T_{i} \ , & \text{otherwise\text{.}}
  \end{cases}$

    \item[3.] $ \overline{v_{i}}\left( \neg_{c_{k}} \ \varphi \right)=\begin{cases}
    \overline{v_{i}}\left( \varphi \right) \ , & \text{If $i$ is not a symbol of $c_{k}$ }\\
    T_{i} \ , & \text{if $i$ is a symbol of $c_{k}$ and $\overline{v_{i}}\left( \varphi \right)=F_{i}$}\\
    F_{i} \ , & \text{if $i$ is a symbol of $c_{k}$ and $\overline{v_{i}}\left( \varphi \right)=T_{i}$\text{.}}
  \end{cases}$

   \item[4.]  $ \overline{v_{i}}\left( \varphi \wedge_{(n)} \psi \right)=\begin{cases}
    T_{i} \ , & \text{if $\overline{ v_{i}}\left( \varphi \right)=\overline{ v_{i}}\left( \psi \right)= T_{i}$}\\
    F_{i} \ , & \text{otherwise\text{.}}
  \end{cases}$

   \item[5.]  $ \overline{v_{i}}\left( \varphi \vee_{(n)} \psi \right)=\begin{cases}
    F_{i} \ , & \text{if $\overline{ v_{i}}\left( \varphi \right)=\overline{ v_{i}}\left( \psi \right)= F_{i}$}\\
    T_{i} \ , & \text{otherwise\text{.}}
  \end{cases}$ 

  \item[6.]  $ \overline{v_{i}}\left( \varphi \longrightarrow_{(n)} \psi \right)=\begin{cases}
    F_{i} \ , & \text{if $\overline{ v_{i}}\left( \varphi \right)= T_{i} $ and $ \overline{ v_{i}}\left( \beta \right)= F_{i}$}\\
    T_{i} \ , & \text{otherwise\text{.}}
  \end{cases}$  

\end{enumerate}

\end{definition}

Each valuation $\overline{v_{i}}$ behaves like the classical valuation in $\textbf{CPC}$ for the connectives and propositional letters, except that we also have the elements $\perp_{c_{k}}$ and the connectives $\neg_{c_{k}}$. The negations behave like classical negation, and $\perp_{c_{k}}$ behaves like classical contradiction. When $c_{k}=(n)$, we have $\overline{v_{i}}(\perp_{(n)})=F_{i}$, and for $\neg_{c_{k}}\varphi$ the value is $V_{i}$ if $\overline{v_{i}}(\varphi)=F_{i}$ and $F_{i}$ if $\overline{v_{i}}(\varphi)=T_{i}$. In definition \ref{def: ivaluation} and in the rest of the paper, all worlds are distinct; the sets in $\{ \mathcal{M}_{i}\}_{1 \leq i \leq n}$ are pairwise disjoint.

\begin{definition}
A valuation $\mathcal{V}$ of $\Omega_{n}$ is a function $\mathcal{V} \colon \Omega_{n} \to \prod_{i=1}^{n}\mathcal{M}_{i}$ such that $\mathcal{V}\left( \varphi \right)= \left( v_{1}\left(\varphi \right), \ldots, v_{n}\left(\varphi \right) \right)$.   
\end{definition}

\begin{definition}
 Each valuation $\mathcal{V}$ extends to $\mathcal{F}\left( \Omega_{n} \right)$ as follows. For each chain $c_{k}$ over $[n]$, define $\overline{\mathcal{V}} \colon \mathcal{F}\left( \Omega_{n} \right) \to \prod_{i=1}^{n} \mathcal{M}_{i}$ by:

\begin{enumerate}

    \item[1.] $\overline{\mathcal{ V}} \left( \varphi \right)= \left( \overline{v_{1}}\left( \varphi \right), \ldots , \overline{v_{n}}\left( \varphi \right)  \right) = \left( v_{1} \left( \varphi \right), \ldots, v_{n}\left( \varphi \right) \right)$ if $\varphi \in P_{n}$; 

    \item[2.]  $\overline{\mathcal{V}} \left( \perp_{c_{k}} \right)= \left( \overline{v_{1}}\left( \perp_{c_{k}} \right), \ldots , \overline{v_{n}}\left( \perp_{c_{k}} \right)  \right)$;
     \item[3.] $\overline{\mathcal{V}} \left( \neg_{c(k)} \varphi \right)= \left( \overline{v_{1}}\left( \neg_{c(k)} \varphi \right), \ldots , \overline{v_{n}}\left(  \neg_{c(k)} \varphi \right)  \right)$;

    \item[4.] $\overline{\mathcal{V}} \left( \varphi \wedge_{(n)} \psi \right)= \left( \overline{v_{1}}\left( \varphi \wedge_{(n)} \psi\right), \ldots , \overline{v_{n}}\left( \varphi \wedge_{(n)} \psi \right)  \right)$; 

    \item[5.] $\overline{\mathcal{V}} \left( \varphi \lor_{(n)} \psi \right)= \left( \overline{v_{1}}\left( \varphi \vee_{(n)}  \psi \right), \ldots , \overline{v_{n}}\left( \varphi \vee_{(n)} \psi \right)  \right)$; 

    \item[6.] $\overline{\mathcal{V}} \left( \varphi \longrightarrow_{(n)} \psi \right)= \left( \overline{v_{1}}\left( \varphi \longrightarrow_{(n)}  \psi \right), \ldots , \overline{v_{n}}\left( \varphi \longrightarrow_{(n)} \psi \right)  \right)$.
   
\end{enumerate}
\end{definition}

\begin{definition}
\label{def: generaltautology}
In $\textup{ \textbf{CPN}}_{n}$, $\varphi$ is a tautology, $\models_{(n)} \varphi$, if and only if for every valuation $\overline{\mathcal{V}}$ we have $\overline{\mathcal{V}}\left( \varphi \right) = \left( T_{1}, \ldots, T_{n} \right)$.
\end{definition}

Definition \ref{def: generaltautology} states that a proposition in $\textbf{CPN}_{n}$ is a tautology if it is a tautology in each of the $n$ worlds. Similarly, we define the analogous notion for the case in which a proposition $\varphi$ takes the value $F_{i}$ in each world.

\begin{definition}
\label{def: generalcontradiction}
 In $\textup{ \textbf{CPN}}_{n}$, $\varphi$ is a contradiction if and only if for every valuation $\overline{\mathcal{V}}$ we have $\overline{\mathcal{V}}\left( \varphi \right) = \left( F_{1}, \ldots, F_{n} \right)$.
\end{definition}

\begin{definition}
Let $\varphi$ be a proposition and let $c_{k}$ be any chain over $[n]$. We say that $\varphi$ is $c_{k}$-\emph{contingent} if $\overline{\mathcal{V}}\left( \varphi \right)$ takes the value $F_{j}$ at all coordinates $j$ such that $j$ is a symbol of $c_{k}$, and the value $T_{i}$ at all coordinates $i$ such that $i$ is a symbol of the complementary chain $c_{n-k}^{'}$.
\end{definition}

After defining the semantics, our goal is to prove the Soundness Theorem and the Completeness Theorem for $\textbf{CPN}_{n}$, which parallel the classical theorems for $\textbf{CPC}$. To this end we first prove the following lemma.

\begin{lemma}
\label{lema: lemasoundness}
 Every axiom is a tautology. 

\begin{proof}
It is easy to check that axioms (A1), (A2), and (A3) are tautologies. We give proofs only for (A4), (A5), (A6), and (A7).
\vspace{0.2cm}

\begin{enumerate}

\item[\textbf{(A4)}]  Suppose there is a valuation $\overline{\mathcal{V}}$ such that for some $i \in [n]$ we have
$$ \overline{v_{i}} \left( \varphi \longrightarrow_{(n)} \left( \perp_{c_{k}} \longrightarrow_{(n)} \neg_{c_{k}}\varphi \right) \right)= F_{i} $$

By definition this means $\overline{v_{i}}(\varphi)=T_{i}$, $\overline{v_{i}}(\perp_{c_{k}})=T_{i}$, and $\overline{v_{i}}(\neg_{c_{k}}\varphi)=F_{i}$. Consider two cases depending on whether $i$ is a symbol of $c_{k}$:
\vspace{0.2cm}
\begin{enumerate}
\item[(a)] If $i$ is a symbol of $c_{k}$, then $\overline{v_{i}}(\perp_{c_{k}})=F_{i}$, a contradiction.
\item[(b)] If $i$ is not a symbol of $c_{k}$, then $\overline{v_{i}}(\neg_{c_{k}}\varphi)=\overline{v_{i}}(\varphi)=T_{i}$, a contradiction.
\end{enumerate}
Hence $\models_{(n)} \varphi \longrightarrow_{(n)} \left( \perp_{c_{k}} \longrightarrow_{(n)} \neg_{c_{k}}\varphi \right)\text{.}$

\item[\textbf{(A5)}] 
Let $c_{k}$ and $c_{r}$ be arbitrary chains over $[n]$. Suppose there is a valuation $\overline{\mathcal{V}}$ and some $i\in[n]$ with

$$ \overline{v_{i}} \left( \neg_{c_{k}}\neg_{c_{r}} \varphi \longleftrightarrow_{(n)} \neg_{c_{k} \otimes c_{r}} \varphi \right)=F_{i} $$

By definition either $\overline{v_{i}}\left( \neg_{c_{k}}\neg_{c_{r}}\varphi \right)=F_{i}$ and $\overline{v_{i}}\left( \neg_{c_{k}\otimes c_{r}}\varphi \right)=T_{i}$, or $\overline{v_{i}}\left( \neg_{c_{k}}\neg_{c_{r}}\varphi \right)=T_{i}$ and $\overline{v_{i}}\left( \neg_{c_{k}\otimes c_{r}}\varphi \right)=F_{i}$. It suffices to consider the first case. We examine four possibilities for $i$ according to whether $i$ is or is not a symbol of $c_{k}$ and $c_{r}$:

\begin{enumerate}
\item[(a)] If $i$ is a symbol of $c_{k}$ and not of $c_{r}$, then $i$ is a symbol of $c_{k}\otimes c_{r}$. From $\overline{v_{i}}\left( \neg_{c_{k}\otimes c_{r}}\varphi \right)= T_{i}$ we get $\overline{v_{i}}\left( \varphi \right)=F_{i}$, while $\overline{v_{i}} \left( \neg_{c_{r}}\varphi \right)= F_{i}$ and hence $\overline{v_{i}}\left( \neg_{c_{k}}\neg_{c_{r}}\varphi \right)= T_{i}$, which contradicts the hypothesis.
\item[(b)] If $i$ is not a symbol of $c_{k}$ and is a symbol of $c_{r}$, then $i$ is a symbol of $c_{k}\otimes c_{r}$. From $\overline{v_{i}}\left( \neg_{c_{k}\otimes c_{r}}\varphi \right)= T_{i}$ we get $\overline{v_{i}}\left( \varphi \right)=F_{i}$, but $\overline{v_{i}} \left( \neg_{c_{r}}\varphi \right)= T_{i}$ and therefore $\overline{v_{i}}\left( \neg_{c_{k}}\neg_{c_{r}}\varphi \right)= T_{i}$, which contradicts the hypothesis.
\item[(c)] If $i$ is a symbol of both $c_{k}$ and $c_{r}$, then $i$ is not a symbol of $c_{k}\otimes c_{r}$. From $\overline{v_{i}}\left( \neg_{c_{k} \otimes c_{r}}\varphi \right)= T_{i}$ we get $\overline{v_{i}} \left( \varphi \right)=T_{i}$, which gives $\overline{v_{i}} \left( \neg_{c_{r}}\varphi \right)= T_{i}$ and hence $\overline{v_{i}}\left( \neg_{c_{k}}\neg_{c_{r}}\varphi \right)= T_{i}$, which contradicts the hypothesis.
\item[(d)] If $i$ is a symbol of neither $c_{k}$ nor $c_{r}$, then $i$ is not a symbol of $c_{k}\otimes c_{r}$. From $\overline{v_{i}} \left( \neg_{c_{k} \otimes c_{r}}\varphi \right)= T_{i}$ we get $\overline{v_{i}}\left( \varphi \right)= T_{i}$, which gives $\overline{v_{i}} \left( \neg_{c_{r}}\varphi \right)= T_{i}$ and hence $\overline{v_{i}}\left( \neg_{c_{k}}\neg_{c_{r}}\varphi \right)= T_{i}$, which contradicts the hypothesis.
\end{enumerate}
Hence $\models_{(n)} \neg_{c_{k}}\neg_{c_{r}} \varphi \longleftrightarrow_{(n)} \neg_{c_{k}\otimes c_{r}}\varphi\text{.}$

\item[\textbf{(A6)}] Let $c_{k}$ and $c_{r}$ be arbitrary chains over $[n]$. Suppose there is a valuation $\overline{\mathcal{V}}$ and some $i\in[n]$ with $$ \overline{v_{i}} \left( \neg_{c_{k}} \perp_{c_{r}} \longleftrightarrow_{(n)} \perp_{c_{k}\otimes c_{r}} \right)=F_{i} $$

By definition either $\overline{v_{i}}\left( \neg_{c_{k}}\perp_{c_{r}} \right)=F_{i}$ and $\overline{v_{i}}\left( \perp_{c_{k}\otimes c_{r}} \right)=T_{i}$, or $\overline{v_{i}}\left( \neg_{c_{k}}\perp_{c_{r}} \right)=T_{i}$ and $\overline{v_{i}}\left( \perp_{c_{k}\otimes c_{r} } \right)=F_{i}$. It suffices to consider the first case. We examine four possibilities for $i$ according to whether $i$ is or is not a symbol of $c_{k}$ and $c_{r}$.
\begin{enumerate}

\item[(a)] If $i$ is a symbol of $c_{k}$ and not a symbol of $c_{r}$, then $i$ is a symbol of $c_{k}\otimes c_{r}$. By definition $\overline{v_{i}}\left( \perp_{c_{k}\otimes c_{r}} \right)= F_{i}$, which contradicts the hypothesis.

\item[(b)] If $i$ is not a symbol of $c_{k}$ and is a symbol of $c_{r}$, then $i$ is a symbol of $c_{k}\otimes c_{r}$. By definition $\overline{v_{i}}\left( \perp_{c_{k} \otimes c_{r}} \right)= F_{i}$, which contradicts the hypothesis.

\item[(c)] If $i$ is a symbol of both $c_{k}$ and $c_{r}$, then $i$ is not a symbol of $c_{k}\otimes c_{r}$. By definition $\overline{v_{i}}\left( \perp_{c_{k}\otimes c_{r}} \right)= T_{i}$, but $\overline{v_{i}} \left( \perp_{c_{r}} \right)=F_{i}$ and $\overline{v_{i}} \left( \neg_{c_{k}}\perp_{c_{r}} \right)=T_{i}$, which contradicts the hypothesis.

\item[(d)] If $i$ is not a symbol of either $c_{k}$ or $c_{r}$, then $i$ is not a symbol of $c_{k}\otimes c_{r}$. By definition $\overline{v_{i}} \left( \perp_{c_{k}\otimes c_{r}} \right)= T_{i}$, then $\overline{v_{i}}\left( \varphi \right)= T_{i}$, but $\overline{v_{i}}\left( \perp_{c_{r}} \right)= T_{i}$ and $\overline{v_{i}}\left( \neg_{c_{k}}\perp_{c_{r}} \right)= T_{i}$, which contradicts the hypothesis.

\end{enumerate}

Hence $\models_{(n)} \neg_{c_{k}} \perp_{c_{r}} \longleftrightarrow_{(n)} \perp_{c_{k}\otimes c_{r}} \text{.}$

\item[\textbf{(A7)}] Let $c_{k}$ be any chain and $c_{r}$ a subchain. Suppose there is a valuation $\overline{\mathcal{V}}$ and some $i\in[n]$ with $$\overline{v_{i}} \left( \perp_{c_{k}} \longrightarrow_{(n)} \perp_{c_{r}} \right)= F_{i}$$

By definition this means $\overline{v_{i}} \left( \perp_{c_{k}} \right)= T_{i}$ and $\overline{v_{i}} \left( \perp_{c_{r}} \right)=F_{i}$. Consider the possibilities according to whether $i$ is or is not a symbol of $c_{k}$ and $c_{r}$.

\begin{enumerate}
   \item[(a)] If $i$ is a symbol of $c_{r}$, then $i$ is also a symbol of $c_{k}$; by definition $\overline{v_{i}}\left( \perp_{c_{k}} \right)= F_{i}$ and $\overline{v_{i}}\left( \perp_{c_{r}} \right)= F_{i}$, which contradicts the hypothesis.

\item[(b)] If $i$ is a symbol of $c_{k}$ but not of $c_{r}$, then $\overline{v_{i}}\left( \perp_{c_{k}} \right)= F_{i}$ and $\overline{v_{i}}\left( \perp_{c_{r}} \right)= T_{i}$, which contradicts the hypothesis.

\item[(c)] If $i$ is not a symbol of $c_{k}$ and not a symbol of $c_{r}$, then $\overline{v_{i}}\left( \perp_{c_{k}} \right)= T_{i}$ and $\overline{v_{i}}\left( \perp_{c_{r}} \right)= T_{i}$, which contradicts the hypothesis.
\end{enumerate}
Hence $\models_{(n)} \perp_{c_{k}} \longrightarrow_{(n)} \perp_{c_{r}} \text{.}$
\end{enumerate}
\end{proof}

\end{lemma}

\begin{theorem}[Soundness]\label{teo: soundness} For every $n\in\mathbb{N}$, if  $\vdash_{(n)} \varphi$ then $\models_{(n)} \varphi$. 

\begin{proof}
 Assume $\vdash_{(n)} \varphi$. We proceed by induction on the length of the derivation of $\varphi$.

\begin{enumerate}
   \item[(a)] If the derivation of $\varphi$ has length 1, then $\varphi$ is an axiom. By Lemma \ref{lema: lemasoundness}, $\varphi$ is a tautology.
    \item[(b)] Suppose the derivation of $\varphi$ has length greater than 1. By the induction hypothesis every formula preceding $\varphi$ in the derivation is a tautology; we show that $\varphi$ is a tautology as well. If $\varphi$ is not an axiom, there is a formula $\psi$ in the derivation such that $\varphi$ is obtained by applying $MP_{n}$: the derivation contains both $\psi \longrightarrow_{(n)} \varphi$ and $\psi$, so ${\psi \longrightarrow_{(n)} \varphi,\ \psi}\vdash_{(n)} \varphi$. By the induction hypothesis, for every valuation $\overline{\mathcal{V}}$ we have $\overline{\mathcal{V}}(\psi)=(T_{1},\ldots,T_{n})$ and $\overline{\mathcal{V}}(\psi \longrightarrow_{(n)} \varphi)=(T_{1},\ldots,T_{n})$. Hence for each $1\le i\le n$ and each $i$-valuation $\overline{v_{i}}$, $\overline{v_{i}}(\psi)=T_{i}$ and $\overline{v_{i}}(\psi \longrightarrow_{(n)} \varphi)=T_{i}$. By the semantics of $\longrightarrow_{(n)}$ it follows that $\overline{v_{i}}(\varphi)=T_{i}$ for all $i$. Therefore $\overline{\mathcal{V}}(\varphi)=(T_{1},\ldots,T_{n})$ for every valuation $\overline{\mathcal{V}}$, so $\models_{(n)} \varphi$.
\end{enumerate} 
\end{proof}    
\end{theorem}
Now that we have proved the soundness theorem, we present some schemes that are not valid in $\textup{\textbf{CPN}}_{n}$, including the law of non-contradiction, the principle of explosion, and the law of excluded middle for weak negations.

\begin{proposition}
\label{prop: firstnoderivable}
For every $n\in\mathbb{N}$, if $c_{k}$ and $c_{r}$ are chains over $[n]$ with $1\leq k,r\leq n-1$, then the following schemes are not derivable in $\textup{\textbf{CPN}}_{n}$:
\begin{enumerate}
    \item[\textup{(i)}] $\varphi \land_{(n)} \neg_{c_{k}}\varphi \longrightarrow_{(n)} \psi $;  
    \item[\textup{(ii)}] $\varphi \land_{(n)} \neg_{c_{k}}\varphi \longrightarrow_{(n)} \perp_{(n)}$; 
    \item[\textup{(iii)}] $\neg_{c_{k}}\left( \varphi \land_{(n)} \neg_{c_{k}}\varphi \right)$; 
    \item[\textup{(iv)}] $\varphi \lor_{(n)} \neg_{c_{k}} \varphi$; 
    \item[\textup{(v)}] $\perp_{c_{k}} \longrightarrow_{(n)} \varphi$;
    \item[\textup{(vi)}] $\neg_{c_{r}} \neg_{c_{k}}\varphi \longrightarrow_{(n)} \varphi$ if $c_{r} \not=c_{k}\text{.}$ 
\end{enumerate}

\begin{proof}
 We prove items (i)–(iv). Items (v) and (vi) are proved similarly.
 \begin{enumerate}
     \item[(i)] Since $c_{k}$ has at most $n-1$ symbols there is $i\in[n]$ such that $i$ is not a symbol of $c_{k}$. Set $\overline{v_{i}}(\psi)=F_{i}$ and $\overline{v_{i}}(\varphi)=T_{i}$. Then $\overline{v_{i}}(\neg_{c_{k}}\varphi)=T_{i}$ and $\overline{v_{i}}(\varphi \land_{(n)} \neg_{c_{k}}\varphi)=T_{i}$. By definition $\overline{v_{i}}(\varphi \land_{(n)} \neg_{c_{k}}\varphi \longrightarrow_{(n)} \psi)=F_{i}$. Hence there is a valuation $\overline{\mathcal{V}}$ with  $$\overline{\mathcal{V}}\left( \varphi \land_{(n)} \neg_{c_{k}}\varphi \longrightarrow_{(n)} \psi \right)= \left(\ldots, F_{i},\ldots \right)$$ By Theorem \ref{teo: soundness} we conclude that $\not\vdash_{(n)} \varphi \land_{(n)} \neg_{c_{k}}\varphi \longrightarrow_{(n)} \psi$. 
      \vspace{0.5cm}
 
     \item[(ii)] Since $c_{k}$ has at most $n-1$ symbols there is $i\in[n]$ such that $i$ is not a symbol of $c_{k}$. Set $\overline{v_{i}}\left( \varphi \right)=T_{i}$. Then $\overline{v_{i}}\left( \neg_{c_{k}}\varphi \right)=T_{i}$ and $\overline{v_{i}}\left( \varphi \land_{(n)} \neg_{c_{k}}\varphi \right)=T_{i}$. Finally $\overline{v_{i}}\left( \perp_{(n)} \right)=F_{i}$ and by definition $\overline{v_{i}}\left( \varphi \land_{(n)} \neg_{c_{k}}\varphi \longrightarrow_{(n)} \perp_{(n)} \right)=F_{i}$. Hence there is a valuation $\overline{\mathcal{V}}$ with $$\overline{\mathcal{V}}\left( \varphi \land_{(n)} \neg_{c_{k}}\varphi \longrightarrow_{(n)} \perp_{(n)} \right)= \left(\ldots,F_{i},\ldots \right)$$ 
    By Theorem \ref{teo: soundness} we conclude that $\not\vdash_{(n)} \varphi \land_{(n)} \neg_{c_{k}}\varphi \longrightarrow_{(n)} \perp_{(n)}$.
     \vspace{0.5cm}

     \item[(iii)] Since $c_{k}$ has at most $n-1$ symbols there is $i\in[n]$ such that $i$ is not a symbol of $c_{k}$. Setting $\overline{v_{i}}\left( \varphi \right)=F_{i}$, by definition $\overline{v_{i}}\left( \varphi \land_{(n)} \neg_{c_{k}}\varphi \right)=F_{i}$. Again by definition $\overline{v_{i}}\left( \neg_{c_{k}} \left( \varphi \land_{(n)} \neg_{c_{k}}\varphi \right) \right)=F_{i}$. Hence there is a valuation $\overline{\mathcal{V}}$ with 

     $$ \overline{\mathcal{V}} \left( \neg_{c_{k}}\left( \varphi \land_{(n)} \neg_{c_{k}} \varphi \right) \right)= \left(\ldots, F_{i}, \ldots \right) $$
     By Theorem \ref{teo: soundness} we conclude that $\not\vdash_{(n)} \neg_{c_{k}}\left( \varphi \land_{(n)} \neg_{c_{k}} \varphi \right) $.
     \vspace{0.5cm}

    \item[(iv)] Since $c_{k}$ has at most $n-1$ symbols there is $i\in[n]$ such that $i$ is not a symbol of $c_{k}$. Set $\overline{v_{i}}\left( \varphi \right)=F_{i}$. By definition $\overline{v_{i}} \left( \neg_{c_{k}}\varphi \right)=F_{i}$ and $\overline{v_{i}} \left( \varphi \lor_{(n)} \neg_{c_{k}}\varphi \right)=F_{i}$. Hence there is a valuation $\overline{\mathcal{V}}$ with $$ \overline{\mathcal{V}} \left( \varphi \lor_{(n)} \neg_{c_{k}}\varphi \right)=\left(\ldots, F_{i}, \ldots \right)$$ By Theorem \ref{teo: soundness} we conclude that $\not\vdash_{(n)} \varphi \lor_{(n)} \neg_{c_{k}}\varphi $.       
 \end{enumerate}
\end{proof}
\end{proposition}

In \cite{paracon} a logical system is defined to be paraconsistent if there exists $\varphi$ such that $\varphi, \neg \varphi \not\vdash \psi$. For weak negations $\neg_{c_{k}}$ where the chain length is strictly less than $n$, we have $\varphi, \neg_{c_{k}}\varphi \not\vdash_{(n)} \psi$, and hence $\textup{\textbf{CPN}}_{n}$ is paraconsistent.

\begin{proposition}
For every $n\in\mathbb{N}$, if $c_{k}$ is a chain over $[n]$ with $1\le k\le n-1$, then the following schemes are not derivable in $\textup{\textbf{CPN}}_{n}:$

\begin{enumerate}
    \item[\textup{(i)}] $\neg_{c_{k}}\varphi \longrightarrow_{(n)} \left( \varphi\longrightarrow_{(n)} \psi \right)$; 
    \item[\textup{(ii)}] $\left( \neg_{c_{k}}\varphi \longrightarrow_{(n)} \neg_{c_{k}}\psi \right) \longrightarrow_{(n)} \left( \psi \longrightarrow_{(n)} \varphi \right)$; 
    \item[\textup{(iii)}] $\left( \varphi \longrightarrow_{(n)} \psi \right) \longrightarrow_{(n)} \left( \neg_{c_{k}}\psi \longrightarrow_{(n)} \neg_{c_{k}}\varphi \right)$; 
    \item[\textup{(iv)}] $\left( \varphi \lor_{(n)} \psi \right) \land_{(n)} \neg_{c_{k}}\varphi \longrightarrow_{(n)} \psi$; 
    \item[\textup{(v)}] $\left( \varphi \longrightarrow_{(n)} \psi \right) \land_{(n)} \neg_{c_{k}}\psi \longrightarrow_{(n)} \neg_{c_{k}}\varphi$; 
    \item[\textup{(vi)}] $\neg_{c_{k}}\varphi \longrightarrow_{(n)} \neg_{c_{k}}\left( \varphi \land_{(n)} \psi \right)$; 
    \item[\textup{(vii)}] $\neg_{c_{k}}\varphi \lor_{(n)} \neg_{c_{k}}\psi \longrightarrow_{(n)} \neg_{c_{k}} \left( \varphi \land_{(n)} \psi \right)$; 
    \item[\textup{(viii)}] $\neg_{c_{k}} \left( \varphi \lor_{(n)} \psi \right) \longrightarrow_{(n)} \neg_{c_{k}}\varphi \land_{(n)} \neg_{c_{k}}\psi$. 
    \end{enumerate}

\begin{proof}
We present only the proofs of items (vii) and (viii). These are the missing implications of De Morgan’s laws for weak negations and are not derivable in $\textbf{CPN}_{n}$.
\begin{enumerate}
    \item[(i)] Since $c_{k}$ has at most $n-1$ symbols there is $i\in[n]$ such that $i$ is not a symbol of $c_{k}$. Set $\overline{v_{i}}\left( \varphi \right)=T_{i}$ and $\overline{v_{i}}\left( \psi \right)=F_{i}$. By definition $\overline{v_{i}}\left( \neg_{c_{k}}\varphi \right)=T_{i}$ and $\overline{v_{i}}\left( \neg_{c_{k}}\psi \right)=F_{i}$. Hence $\overline{v_{i}}\left( \neg_{c_{k}}\varphi \lor_{(n)} \neg_{c_{k}}\psi \right)=T_{i}$, $\overline{v_{i}}\left( \varphi \land_{(n)} \psi \right)=F_{i}$, and $\overline{v_{i}}\left( \neg_{c_{k}}\left( \varphi \land_{(n)} \psi \right) \right)=F_{i}$. Then
    $$ \overline{v_{i}}\left( \neg_{c_{k}}\varphi \lor_{(n)} \neg_{c_{k}}\psi \longrightarrow_{(n)} \neg_{c_{k}}\left( \varphi \land_{(n)} \psi \right) \right)= F_{i} $$
   This means there is a valuation $\overline{\mathcal{V}}$ with
    $$\overline{\mathcal{V}}\left( \neg_{c_{k}} \varphi \lor_{(n)} \neg_{c_{k}}\psi \longrightarrow_{(n)} \neg_{c_{k}} \left( \varphi \land_{(n)} \psi \right)\right)= \left( \ldots, F_{i},\ldots\right)$$
    By Theorem \ref{teo: soundness} we conclude that $\not\vdash_{(n)} \neg_{c_{k}} \varphi \lor_{(n)} \neg_{c_{k}}\psi \longrightarrow_{(n)} \neg_{c_{k}} \left( \varphi \land_{(n)} \psi \right)$.
    \vspace{0.5cm}

    \item[(ii)]Since $c_{k}$ has at most $n-1$ symbols there is $i\in[n]$ such that $i$ is not a symbol of $c_{k}$. Set $\overline{v_{i}}\left( \varphi \right)=F_{i}$ and $\overline{v_{i}}\left( \psi \right)=T_{i}$. By definition $\overline{v_{i}}\left( \neg_{c_{k}}\varphi \right)=F_{i}$ and $\overline{v_{i}}\left( \neg_{c_{k}}\psi \right)=T_{i}$. Hence $\overline{v_{i}}\left( \neg_{c_{k}}\varphi \land_{(n)} \neg_{c_{k}}\psi \right)=F_{i}$, $\overline{v_{i}}\left( \varphi \lor_{(n)} \psi \right)=T_{i}$, and $\overline{v_{i}}\left( \neg_{c_{k}}\left( \varphi \lor_{(n)} \psi \right) \right)=T_{i}$. Then
    $$ \overline{v_{i}}\left( \neg_{c_{k}}\left( \varphi \lor_{(n)} \psi \right)  \longrightarrow_{(n)} \neg_{c_{k}}\varphi \land_{(n)} \neg_{c_{k}}\psi \right)= F_{i} $$
    This means there is a valuation $\overline{\mathcal{V}}$ with
    $$\overline{\mathcal{V}}\left( \neg_{c_{k}}\left( \varphi \lor_{(n)} \psi \right)  \longrightarrow_{(n)} \neg_{c_{k}}\varphi \land_{(n)} \neg_{c_{k}}\psi\right)= \left( \ldots, F_{i},\ldots\right)$$
    By Theorem \ref{teo: soundness} we conclude that $\not\vdash_{(n)} \neg_{c_{k}}\left( \varphi \lor_{(n)} \psi \right)  \longrightarrow_{(n)} \neg_{c_{k}}\varphi \land_{(n)} \neg_{c_{k}}\psi$. 
\end{enumerate}
\end{proof}
\end{proposition}
\section{Completeness for $\textbf{CPN}_{n}$}

In this final section we present the Completeness Theorem for $\textbf{CPN}_{n}$, whose proof is similar to that in $\textbf{CPC}$. For further references see \cite{caicedo}, \cite{Hod}, and \cite{Men}.

\begin{definition}
Let $\mathcal{V}$ be a fixed valuation. For each formula $\varphi$ and each chain $c_{k}$ over $[n]$, define:
 $$\varphi^{\mathcal{V}}:= \neg_{c_{k}}\varphi \hspace{0.5cm} \text{if} \ \varphi \ \text{is} \ c_{k}-\text{contingent}\text{.} $$
\end{definition}

\begin{example}
In $\textbf{CPN}_{2}$, for any formula $\varphi$:
$$ \varphi^{\mathcal{V}}=\begin{cases}
   \varphi , & \text{if $\overline{\mathcal{V}}\left( \varphi \right)= \left(T_{1}, T_{2} \right)$ \text{.}}\\
    \neg_{1}\varphi \ , & \text{if $\overline{\mathcal{V}}\left( \varphi \right)= \left(F_{1}, T_{2} \right)$ \text{.}}\\
    \neg_{2}\varphi \ , & \text{if $\overline{\mathcal{V}}\left( \varphi \right)= \left(T_{1}, F_{2} \right)$ \text{.}} \\
    \neg_{(2)}\varphi \ , & \text{if $\overline{\mathcal{V}}\left( \varphi \right)= \left(F_{1}, F_{2} \right)$ \text{.}}
  \end{cases}$$
\end{example}

\begin{lemma}
 For every $n\in\mathbb{N}$ and every chain $c_{s}$ over $[n]$, we have:
\begin{enumerate}
\item[\textup{(i)}] $\varphi^{\mathcal{V}} \vdash_{(n)} \left( \neg_{c_{s}}\varphi \right)^{\mathcal{V}}$;
\item[\textup{(ii)}] $\varphi^{\mathcal{V}}, \psi^{\mathcal{V}} \vdash_{(n)} \left( \varphi \longrightarrow_{(n)} \psi \right)^{\mathcal{V}}\text{.}$
\end{enumerate}

\begin{proof} For $n\in\mathbb{N}$ we have:
\begin{enumerate}

 \item[(i)]  If $\varphi$ is $c_{k}$-contingent for some chain $c_{k}$ over $[n]$, then $\varphi^{\mathcal{V}}=\neg_{c_{k}}\varphi$. By Proposition \ref{prop: cknegations} and $TD_{n}$ it follows that for any $c_{s}$ over $[n]$, $\neg_{c_{k}}\varphi \vdash_{(n)} \neg_{c_{s}}\neg_{c_{s}}\neg_{c_{k}}\varphi$, or equivalently by axiom \textbf{(A5)}, $\neg_{c_{k}}\varphi \vdash_{(n)} \neg_{c_{k}\otimes c_{s}}\neg_{c_{s}}\varphi$. Hence $\varphi^{\mathcal{V}} \vdash_{(n)} \left( \neg_{c_{s}}\varphi \right)^{\mathcal{V}}\text{.}$

  \item[(ii)] If $\varphi$ is $c_{k}$-contingent and $\psi$ is $c_{s}$-contingent for chains $c_{k}$ and $c_{s}$ over $[n]$, then by Theorem \ref{teo: fundamental} we have $\neg_{c_{k}}\varphi, \neg_{c_{s}}\psi \vdash_{(n)} \neg_{c_{s}\otimes c_{k}} \left( \varphi \longrightarrow_{(n)} \psi \right)$. Hence $\varphi^{\mathcal{V}}, \psi^{\mathcal{V}} \vdash_{(n)} \left( \varphi \longrightarrow_{(n)} \psi \right)^{\mathcal{V}}\text{.}$
 \end{enumerate}
\end{proof}

\end{lemma}

\begin{lemma}
\label{lema: propositionalletters}
 For every $n\in\mathbb{N}$, if $p_{1},p_{2},\ldots,p_{m}$ are the propositional letters occurring in $\varphi$ and $\mathcal{V}$ is a valuation, then $p_{1}^{\mathcal{V}},p_{2}^{\mathcal{V}},\ldots,p_{m}^{\mathcal{V}} \vdash_{(n)} \varphi^{\mathcal{V}}\text{.}$
\begin{proof}
The proof proceeds by induction on the complexity of $\varphi$.
\end{proof} 
\end{lemma}

\begin{theorem}[Completeness] \label{teo: completeness}
For every, if $\models_{(n)} \varphi$ then $\vdash_{(n)} \varphi\text{.}$ 

\begin{proof}
 If $\models_{(n)} \varphi$ then $\varphi$ is a tautology, that is, for every valuation $\mathcal{V}$ we have $\overline{\mathcal{V}}(\varphi)=(T_{1},\ldots,T_{n})$. By definition $\varphi^{\mathcal{V}}=\varphi$. By Lemma \ref{lema: propositionalletters}, if $p_{1},\ldots,p_{m}$ are the propositional letters occurring in $\varphi$, then $p_{1}^{\mathcal{V}},\ldots,p_{m}^{\mathcal{V}} \vdash_{(n)} \varphi$. Now let $\mathcal{V}$ be a valuation on $p_{2},\ldots,p_{m}$; extend it to $p_{1}$ as follows. Since $p_{1}$ has $2^{n}$ truth values, for any chain $c_{k}$ with $0\le k\le n$, if $p_{1}$ is $c_{k}$-contingent then $p_{1}^{\mathcal{V}}=\neg_{c_{k}}p_{1}$. There is also the case that $p_{1}$ is $c_{n-k}^{\prime}$-contingent, that is, $p_{1}^{\mathcal{V}}=\neg_{c_{n-k}^{\prime}}p_{1}$, which is equivalent to $p_{1}^{\mathcal{V}}=\neg_{(n)}\neg_{c_{k}}p_{1}$. Hence for every chain $c_{k}$ over $[n]$ with $0\le k\le n$

 $$\neg_{c_{k}}p_{1},p_{2}^{\mathcal{V}},\ldots, p_{m}^{\mathcal{V}} \vdash_{(n)} \varphi $$ and $$ \neg_{(n)}\neg_{c_{k}}p_{1}, p_{2}^{\mathcal{V}},\ldots,p_{n}^{\mathcal{V}} \vdash_{(n)} \varphi $$ 
 By $TD_{n}$ we have $p_{2}^{\mathcal{V}},\ldots,p_{m}^{\mathcal{V}} \vdash_{(n)} \neg_{c_{k}}p_{1} \longrightarrow_{(n)} \varphi$ and $p_{2}^{\mathcal{V}},\ldots,p_{m}^{\mathcal{V}} \vdash_{(n)} \neg_{(n)}\neg_{c_{k}}p_{1} \longrightarrow_{(n)} \varphi$. Moreover, by the classical rule in $\textbf{CPN}_{n}$, which states that from $\psi \longrightarrow_{(n)} \chi$ and $\neg_{(n)}\psi \longrightarrow_{(n)} \chi$ we may infer $\chi$, we obtain 

 $$ p_{2}^{\mathcal{V}},\ldots, p_{n}^{\mathcal{V}} \vdash_{(n)} \varphi$$
After repeating this procedure $m$ times we conclude that $\vdash_{(n)}\varphi\text{.}$
\end{proof}

\end{theorem}

\begin{theorem}
 For every $n\in\mathbb{N}$, if $\vdash_{(n+1)} \varphi$ then $\vdash_{(n)} \varphi$. 

\begin{proof}
 Assume $\not\vdash_{(n)} \varphi$. By Theorem \ref{teo: completeness} we have $\not\models_{(n)} \varphi$. By definition there is a valuation $\overline{\mathcal{V}}$ and some $i\in [n]$ with $\overline{v_{i}}(\varphi)=F_{i}$. Since $[n]\subseteq [n+1]$, we also have $i\in [n+1]$. Again by definition, $\not\models_{(n+1)}\varphi$. By Theorem \ref{teo: soundness} we conclude that $\not\vdash_{(n+1)}\varphi$, a contradiction.   
\end{proof}

\end{theorem}

\renewcommand{\refname}{References}


\end{document}